\documentclass[12pt]{article} 
\input{epsf.sty}
\usepackage{graphicx} 
\usepackage{amsmath}
\usepackage{amssymb} 
\usepackage{amscd}
\usepackage{amsthm} 
\usepackage{epsfig} 
\usepackage{latexsym}
\usepackage{enumerate}

\newtheorem{theorem}{Theorem}[section]
\newtheorem{lemma}[theorem]{Lemma}
\newtheorem{corollary}[theorem]{Corollary} 
 
\newtheorem{proposition}[theorem]{Proposition}

\newtheorem{addendum}[theorem]{Addendum}

\title{A Combination Theorem For Convex Hyperbolic Manifolds, With Applications To Surfaces In 3-Manifolds.}
\author{Mark Baker and Daryl Cooper}
\date{}

\def\gap{\vspace{.3cm}\noindent}

\def\demo{ {\bf Proof.} }
\def\nb{\newline$\bullet\ \ \ $}
\def\smallskip{\vspace{.15cm}}
\def\medskip{\vspace{.3cm}}
\def\text{\mbox}
\def\qed{\hfill$\square$}

\begin{document}
\maketitle
\abstract{We prove the convex combination theorem for hyperbolic $n$-manifolds. Applications are given both in high dimensions and in $3$ dimensions. One consequence is that given two geometrically finite subgroups of a discrete group of isometries of hyperbolic $n$-space, satisfying a natural condition on their parabolic subgroups, there are finite index subgroups which generate a subgroup that is an amalgamated free product. Constructions of infinite volume hyperbolic $n$-manifolds are described by gluing lower dimensional manifolds. It is shown that every slope on a cusp of a hyperbolic $3$-manifold is a multiple immersed boundary slope. If a $3$-manifold contains a maximal surface group not carried by an embedded surface then it contains the fundamental group of a book of $I$-bundles with more than two pages.}
\footnote{AMS subj. class: 57M50 (Primary); 30F40 (Secondary)}

\section{Introduction}
The Klein-Maskit combination theorems \cite{MAS} assert that under certain circumstances two Kleinian groups, $\Gamma_1,\Gamma_2$ corresponding to hyperbolic manifolds $M_1,M_2$ are subgroups of another Kleinian group $\Gamma$ corresponding to a hyperbolic manifold $M$ diffeomorphic to that obtained by gluing $M_1$ and $M_2$ along a boundary component. There is also an HNN version, and an orbifold version.

We prove the {\em convex combination theorem} which allows the gluing of two convex hyperbolic $n$-manifolds along isometric submanifolds which are not necessarily boundary parallel. 
One consequence is the {\em virtual amalgam theorem} which states that that if $G$ and $H$ are two geometrically-finite subgroups of a discrete group in $Isom({\mathbb H}^n),$ and if their parabolic subgroups are {\em compatible} (defined in section 5), then there are finite index subgroups $G',H'$ which generate a geometrically finite subgroup that is the free product of $G'$ and $H'$ amalgamated along $G'\cap H'.$ 

One application is proving the existence of certain kinds of surface groups in hyperbolic 3-manifolds. For example we show that given any slope on a torus boundary component of a compact 3-manifold, $M,$ with hyperbolic interior, there is a compact, immersed, oriented, geometrically-finite, $\pi_1$-injective, surface not homotopic into the boundary of $M$ whose boundary consists of two components each of which wraps the same number of times, but in opposite directions, around the given slope.  In particular $\infty$ is a multiple immersed boundary slope of every hyperbolic knot. 
The first author introduced this concept in \cite{B} and gave an example of a hyperbolic once-punctured torus bundle  with infinitely many immersed boundary slopes.  Oertel \cite{OE} gave the first example of a hyperbolic manifold such that every slope is an immersed boundary slope, and Maher \cite{MA} proved the latter holds for hyperbolic 2-bridge knots and certain other cases.

Another application is that if the fundamental group of a non-Haken closed hyperbolic 3-manifold contains a surface group then it contains the fundamental group of an irreducible, boundary-irreducible, compact 3-mani\-fold (a book of $I$-bundles) with arbitrarily large second betti number.

The basic idea is to glue two convex hyperbolic manifolds together to obtain a hyperbolic manifold which is not in general convex. It is well known that every $3$-manifold that is not closed admits a non-convex hyperbolic metric, thus such metrics are in general too abundant to provide useful information. What is required is a condition which ensures that the result of the gluing can be thickened to be convex. If this can be done then one obtains a discrete subgroup of isometries of hyperbolic space.

The  {\em convex combination theorem}  asserts (roughly speaking) that there is a universal constant $\kappa,$ independent of dimension, such that if two convex hyperbolic manifolds can be glued together in a way that  is compatible with gluing their  $\kappa$-thickenings, then the resulting manifold can be thickened to be convex.

To deduce the virtual amalgam theorem from the convex combination theorem  involves two hyperbolic manifolds, $A$ and $B,$ which are isometrically immersed into a hyperbolic manifold $M.$ One wishes to glue the basepoints of $A$ and $B$ together and the requirement that the identification space is a manifold forces further identifications to be made between $A$ and $B.$ Subgroup separability arguments are used  to ensure that certain finite covers of $A$ and $B$ can be glued so that they embed in the resulting identification space.  In order to satisfy the thickening hypotheses of the convex combination theorem one might need to take large finite covers of the manifolds concerned.  

The paper is organized as follows. In section 2 we discuss  convex hyperbolic manifolds and prove the {\em convex combination theorem} (\ref{convexcombination}). In section 3 we study cusps of geometrically-finite, convex hyperbolic manifolds. In section 4 we introduce the notion of {\em induced gluing} alluded to in the preceding paragraph and prove the {\em virtual simple gluing theorem} (\ref{virtualsimpleglue}) which ensures the manifolds that are glued embed in the resulting space.  In section 5 we prove the {\em virtual amalgam theorem} (\ref{virtualamalgam}) and the {\em virtual convex combination theorem} (\ref{virtualcombination}). In section 6 we construct some higher dimensional convex hyperbolic manifolds of infinite volume by gluing lower dimensional ones. In section 7 we show that certain groups are LERF and extend an argument of Scott's that finitely-generated subgroups of surface groups are almost geometric to the case of finitely-generated separable subgroups of  three-manifold groups.  In section 8 we prove some new results about surface groups in hyperbolic 3-manifolds as well as sketching new proofs of some results of the second author and Long about virtually-Haken Dehn-filling and the existence of surface groups in most Dehn-fillings.  In section 9 we apply these tools to the study of immersed boundary slopes. 

The convex combination theorem is related to work of Bestvina-Feighn \cite{BF}, Gitik \cite{G2} and Dahmani \cite{FD} who proved various combination theorems for (relatively) word hyperbolic groups. The convex combination theorem implies that  certain groups are discrete groups of isometries of hyperbolic space, a conclusion which does not follow from the group-theory results mentioned. It seems possible that there is a common generalization of the Klein-Maskit theorem and the convex combination theorem, but this will have to await a mythical future  paper.

The train of ideas in this paper originated with the work of B. Freedman and M.H. Freedman \cite{FF} who constructed certain closed surfaces by a tubing operation. If the surfaces involved are {\it far enough apart}  in a certain sense (if the tube is  {\it long enough}) then the resulting surface is $\pi_1$-injective. There are by now several proofs of this  and related facts, and this paper provides new ones. 

Both authors thank the Universit$\acute{e}$ de Rennes 1 and UCSB for hospitality and partial support during the time this paper was written. This work was also partially supported by NSF grants DMS0104039 and DMS0405963.

\section{The Convex Combination Theorem.}

In this section we prove the convex combination theorem. This requires a brief discussion of non-convex hyperbolic manifolds. To do this we need to extend some ideas from the more well-known context of convex to that of non-convex hyperbolic manifolds.

\medskip
\noindent{\bf Definition.} A {\em hyperbolic manifold} is a connected manifold with boundary (possibly empty) equipped with a Riemannian metric which is hyperbolic i.e. constant sectional curvature $-1.$ 

\medskip\noindent{\bf Warning:} We do not assume the holonomy of a  hyperbolic manifold is a discrete group of isometries of hyperbolic space. 

\medskip
  We will primarily be interested in the case that the boundary is piecewise smooth, but has corners. 
 Let $\tilde{M}$ denote the universal cover of a hyperbolic manifold  $M.$ There is a local isometry called the {\em developing map}  and a homomorphism of groups called the {\em holonomy}
$$dev:\tilde{M}\rightarrow{\mathbb H}^n\qquad\qquad hol:\pi_1(M)\longrightarrow Isom({\mathbb H}^n)$$
such that for all $x\in\tilde{M}$ and all $g\in\pi_1M$ we have $dev(g\cdot x)=hol(g)\cdot dev(x).$
    
\medskip
\noindent{\bf Definition.}  A hyperbolic manifold  is {\em convex} if every two points in the universal cover  are connected by a geodesic. 

\medskip
In particular the quotient of hyperbolic space by a discrete torsion-free group of isometries is convex. 
The following is  easy to check:

\begin{proposition}[characterize convex]\label{characterizeconvex}\  \\ Suppose that $M$ is a hyperbolic manifold. Then the following are equivalent. \\
(a) $M$ is convex.\\
(b)  Every path in $M$  is homotopic rel endpoints to a geodesic in $M.$ \\
(c) The developing map is injective with image a convex subset of hyperbolic space.
\end{proposition}

\begin{proposition}[convex has injective holonomy]\label{convexinjhol}\  \\ If $M$ is a convex hyperbolic manifold then the holonomy is  injective.\\ Futhermore $M=dev(\tilde{M})/hol(\pi_1M).$\end{proposition}
\demo It follows from (\ref{characterizeconvex})(c)  that the developing map is an isometry onto its image. Since $\pi_1M$ acts freely by isometries on its universal cover,  the holonomy is  injective and  $M$ is isometric to $dev(\tilde{M})/hol(\pi_1M).$ \qed

\begin{proposition}[local isometry from convex is $\pi_1$-injective]\label{convexinjective}
\   \\
Suppose $M$ and $N$ are hyperbolic manifolds, $M$ is convex, and $f:M\rightarrow N$ is a local isometry. Then $f_*:\pi_1M\rightarrow\pi_1N$ is injective. In particular, if $N={\mathbb H}^n$ then $M$ is simply connected.\end{proposition}
\demo It is easy to check that  $hol_M = hol_N\circ f_*.$ Since $M$ is convex,  $hol_M$ is injective by (\ref{convexinjhol}). Thus $f_*$ is injective.\qed

\medskip
Consider two closed geodesics in a hyperbolic surface which intersect in two points. Let $A$ and $B$ be small convex neighborhoods of these geodesic. Then $A\cap B$ is the disjoint union of two discs each of which is convex. More generally we have:

 \begin{lemma}[intersection of closed convex is convex union]\label{convexintersect}\  \\ Suppose $\{ M_i : i\in I \}$  are convex hyperbolic manifolds which are closed subsets of a hyperbolic $n$-manifold $M.$ Then every component of $\cap_{i\in I} M_i$ is a convex hyperbolic manifold.
 \end{lemma}
 \demo Let $p:\tilde{M}\rightarrow M$ be the universal cover and $\tilde{M}_i$ a component of $p^{-1}(M_i).$  Since $M_i$ is convex it follows from (\ref{convexinjective}) applied to $M_i\hookrightarrow M$ that $p|:\tilde{M_i}\rightarrow M_i$ is the universal cover. 
 
 Let $C$ be a component of  $\cap_{i\in I} M_i.$ Let $\tilde{C}$ be a component of $p^{-1}(C).$ 
  Consider the components $\tilde{M}_i\subset p^{-1}M_i$ which contain $\tilde{C}.$
  Since $dev_M$ embeds each $\tilde{M}_i$ it also embeds $K = \cap_i \tilde{M}_i$ into ${\mathbb H}^n.$ Thus $dev_M(K) =  \cap_i\ dev_M(\tilde{M}_i)$ is a closed convex subset of ${\mathbb H}^n$ and therefore a manifold. Hence $K$ is a convex manifold. Clearly $\tilde{C}\subset K.$ Choose $x\in\tilde{C}$ and $y\in K$ then, because $K$ is convex, there is a unique geodesic segment $[x,y]$ in $K$ with endpoints $x$ and $y.$ This segment is in every $\tilde{M}_i$ and therefore $p([x,y])$ is contained in every $M_i$ and thus in $C.$ It follows that $[x,y]$ is contained in $\tilde{C}$ hence $y\in \tilde{C}.$ Thus $\tilde{C}=K.$
 It follows that $C$ is a convex hyperbolic manifold. \qed

\medskip
{\bf Definition.} Suppose $M$ is a convex hyperbolic $n$-manifold and $A$ is a non-empty, connected subset of $M.$ The {\em convex hull}, $CH(A),$ of $A$ is defined to be  the intersection of all convex manifolds in $M$ which are closed subsets of $M$ and which contain $A.$ 

\begin{proposition}[convex hulls are convex]\label{convexhull}\  \\ If $M$ is a convex hyperbolic $n$-manifold and $A$ is a non-empty connected subset of $M$ then every component of $CH(A)$ is a convex manifold of some dimension $k\le n.$\end{proposition}
\demo  Since $A$ is connected there is a unique component, $C,$ of $CH(A)$ which contains $A.$ By (\ref{convexintersect})  $C$ is a convex manifold which contains $A$ so $CH(A)=C.$ \qed

\medskip
 There are many examples of non-convex hyperbolic manifolds. For example an immersion of a punctured torus into the hyperbolic plane induces a pull-back metric on the punctured torus with trivial holonomy.  Every non-compact 3-manifold can be immersed into Euclidean space and hence into ${\mathbb H}^3.$ It follows that every such manifold has a hyperbolic metric.
  
 \medskip
 We want to know when a non-convex manifold corresponds to a discrete group of isometries. The preceding examples  do not have this property. There are several equivalent ways to describe the desired property, and one involves the notion of thickening.
 
 \medskip{\bf Definition.} A hyperbolic $n$-manifold $N$ is a {\em thickening} of a hyperbolic $n$-manifold, $M,$  if $M\subset N$ and $incl_*:\pi_1M\rightarrow\pi_1N$ is an isomorphism.  If, in addition, $N$ is convex then we say $N$ is a {\em convex thickening} of $M.$ 
 The following is easy to check:

\begin{proposition}[convex thickenings]\label{injectivedeveloping}\ \\ Suppose that $M$ is a  hyperbolic $n$-manifold.  Then the following are equivalent.\\
(a) The developing map $dev:\tilde{M}\rightarrow{\mathbb H}^n$ is injective.\\
(b)  The holonomy of $M$ is a discrete torsion-free group $\Gamma\subset Isom({\mathbb H}^n)$ and there is an isometric embedding $f:M\rightarrow N={\mathbb H}^n/\Gamma$ such that $f_*:\pi_1M\rightarrow\pi_1N$ is an isomorphism.\\
(c) $M$ has a convex thickening.\end{proposition}

 \medskip
 We will often use the developing map to identify the universal cover
$\tilde{M}$ of a convex manifold with $dev(\tilde{M})\subset{\mathbb H}^n.$ If
$M$ is a convex hyperbolic manifold and
$K\ge0$ the {\em $K$-thickening} of $M$ is 
$$T_K(M)\ =\ N_K(\tilde{M})/\pi_1M$$ where $N_K(\tilde{M}) = \{x\in{\mathbb H}^n : d(x,\tilde{M})\le K
\}$ is the {\it $K$-neighborhood} of $\tilde{M}$ in ${\mathbb H}^n.$ With this notation $T_{\infty}(M)={\mathbb H}^n/hol(\pi_1M)$ is the geodesically-complete manifold that is a thickening of $M.$  The following is immediate:
 
 \begin{proposition}[thickening]\label{thickening} Suppose $M$ is a convex hyperbolic manifold. Then:\\
 (a) $T_K(M)$ is a convex hyperbolic manifold which contains an
 isometric copy of $M$ and is unique up to isometry fixing $M.$\\(b)  If $x\in M$ and $y\in\partial T_K(M)$ then $d(x,y)\ge K.$\end{proposition}

\medskip
 
The convex combination theorem (\ref{convexcombination}) gives a sufficient condition to ensure that the union of two convex hyperbolic $n$-manifolds has a convex thickening, and so has holonomy a discrete subgroup of $Isom({\mathbb H}^n).$ The following example shows that some additional hypothesis is needed to ensure this.

\medskip
{\bf Example.} Suppose that $S(\ell,\theta,K)=M_1\cup M_2$
is homeomorphic to a punctured torus with an incomplete hyperbolic metric
and that $M_1$ and $M_2$ are hyperbolic annuli isometric to
$K$-neighborhoods of closed geodesics
$\gamma_1,\gamma_2$ of length $\ell.$ Suppose $M_1\cap M_2$ is a
disc and that the angle between the geodesics $\gamma_1,\gamma_2$ is
$\theta.$ Then given $\theta\in(0,\pi),$ the set of
$\ell>0$ for which there is $K>0$ such that the developing map
$dev:\tilde{S}\rightarrow{\mathbb H}^2$ is injective is an interval
$[\ell_0(\theta),\infty)$. By Margulis's theorem there is $\mu>0$ such that
$\ell_0>\mu$ for all $\theta.$ Also $\ell_0(\theta)\to\infty$ as
$\theta\to0.$ On the other hand it
is easy to convince oneself that if $K>100$ then the developing map
is always injective.

 \medskip
 The next result holds even when the manifold $M$ has no convex thickening.

 \begin{corollary}[union of convex gives amalgamated free product]\label{freeproduct} \  \\ Suppose $M=M_1\cup M_2$ is a connected hyperbolic $n$-manifold which is the union of two convex hyperbolic $n$-manifolds $M_1,M_2$ and suppose that $M_1\cap M_2$ is connected. Given a basepoint $x\in M_1\cap M_2$ then\\ $\pi_1(M,x) = \pi_1(M_1,x)*_G\pi_1(M_2,x)$ where $G = \pi_1(M_1\cap M_2,x).$ 
 \end{corollary}
 \demo By (\ref{convexintersect}) $M_1\cap M_2$ is convex. It follows from (\ref{convexinjective}) that $M_1\cap M_2$ is $\pi_1$-injective in $M.$ The result follows from Van Kampen's theorem.\qed
 
\medskip
The following theorem asserts, very roughly,  that there is a universal constant, $\kappa,$ such that if $M$ is a (probably  non-convex) hyperbolic $n$-manifold which is the union of two convex hyperbolic submanifolds and if $M$ has a $\kappa$-thickening with the same topology, then $M$ has a convex thickening. The intuition for this result is lemma (\ref{convexsetlemma}) which says that if two convex sets in hyperbolic space intersect then their convex hull is within a small distance of their union. The fact one can thicken the submanifolds without bumping means that the universal cover is made of convex sets (covers of the submanifolds) which are far apart in some sense. Then the convex hull construction more or less only notices two of the convex sets  at any one time and so the convex hull is close to the union.

\centerline{\epsfysize=60mm
\epsfbox{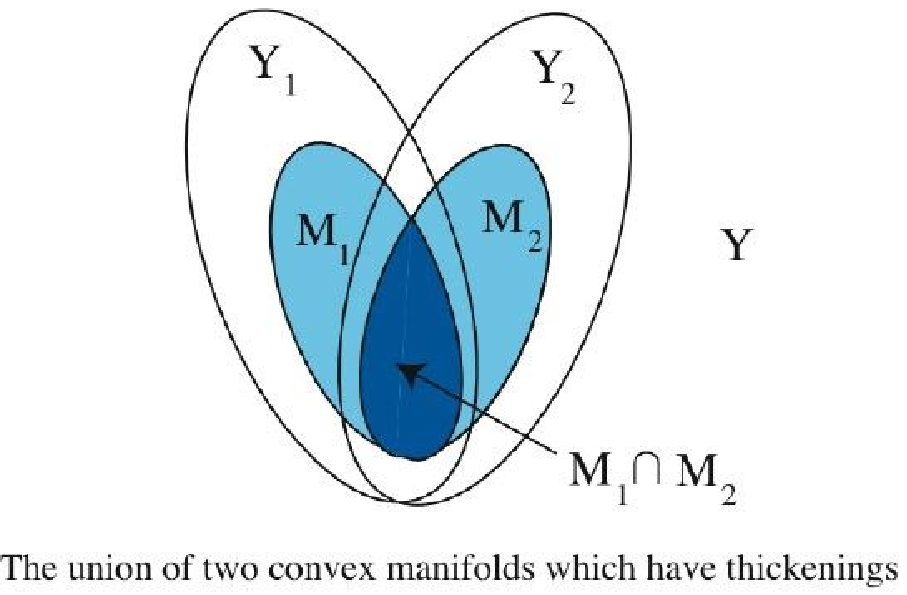}}   

\begin{theorem}[convex combination theorem]\label{convexcombination} \ \\ There is a universal constant, $\kappa,$ called the  {\bf thickening constant}  with the following property. Suppose the following conditions are satisfied:\\
(1)  $Y=Y_1\cup Y_2$ is a connected hyperbolic $n$-manifold which is the union of two convex $n$-submanifolds $Y_1$ and $Y_2.$\\
(2)  $M=M_1\cup M_2$ is a connected hyperbolic $n$-manifold which is the union of two convex $n$-submanifolds $M_1$ and $M_2.$\\
(3)  $M\subset Y$ and $Y_i$ is a thickening of $M_i$ for $i=1,2.$\\
(4) For  all $p\in M$ the exponential map $exp_p:T_p^{\kappa}M\rightarrow Y$ is defined,  where $T^{\kappa}_pM$ denotes the set of tangent vectors in $T_pM$ of length at most $\kappa.$\\
(5) {\bf No bumping:} for $i\in\{1,2\}$  and for all $p\in \overline{M}_i\setminus int(M_1\cap M_2)$ we have $exp_p(T_p^{\kappa}M_i)\subset Y_i.$ \\
(6) Every component of $Y_1\cap Y_2$ contains a point of $M_1\cap M_2.$

\medskip
 \noindent  Then $M$ has a convex thickening;  in other words  there is a hyperbolic $n$-manifold $N={\mathbb H}^n/\pi_1M$ which contains $M$ and $incl_*:\pi_1M\longrightarrow\pi_1N$ is an isomorphism. Furthermore, if $Y$ has finite volume then $M$ is geometrically finite. Also $\kappa\le 6.$ \end{theorem}
\demo   Here is a sketch of the proof that $\kappa=6$ satisfies the theorem;  the details follow. We need to show that the developing map sends the universal cover, $\tilde{M},$ of $M$ injectively into hyperbolic space. To do this we will show that between any two points in $\tilde{M}$ there is a  geodesic in the universal cover, $\tilde{Y},$ of $Y$ connecting them. The image of this geodesic under the developing map is then a geodesic in hyperbolic space and therefore the endpoints are distinct. 

We need to understand the universal cover of $Y$ and its image under the developing map. Convexity of $Y_1$ and $Y_2$  implies that $\tilde{Y}$ is a union of copies of the universal covers of $Y_1$ and $Y_2.$  We show $\pi_1M\cong\pi_1Y,$ and it follows that $\tilde{Y}$ contains the universal cover, $\tilde{M},$ of $M$  and each copy of the universal cover of $Y_i$ contains exactly one copy of the universal cover of $M_i.$ 
 
 To show there is such a geodesic we take a shortest path, $\gamma,$ in the  $2$-neighborhood (contained in $\tilde{Y}$) of $\tilde{M}$ between the two points. Since $M_i$ is convex, near any point of $\gamma$ which is in the interior of the $2$-neighborhood of either $\tilde{M}_1$ or $\tilde{M}_2$ the path is a geodesic. Thus $\gamma$ can only fail to be a geodesic at corners that are on the intersection of the boundaries of the $2$-neighborhoods. So it suffices to show $\gamma$ has no corners.
 
 Condition (5) (no bumping) is used to show that the distance between corners is large (bigger than $\kappa-4$). Thus $\gamma$ is a union of long geodesic segments each of  which starts and ends within a distance of $2$ of some convex set, $\tilde{M}_i.$ It follows that the midpoint of such a segment is then very close (less than $1$) to $\tilde{M}_i.$ We deduce that near a corner of $\gamma$ there is a subpath, $\gamma',$ which consists of two long geodesic segments that meet at that corner and one  endpoint of $\gamma'$ is very close to some $\tilde{M}_1$ and the other to some $\tilde{M}_2.$ 
 
The  convex sets $\tilde{M}_1,\tilde{M}_2$ intersect near the corner. In hyperbolic space the union of two convex sets which intersect is nearly convex:  every point in the convex hull of the union lies within a distance $1$ of the union of the convex sets. Thus the geodesic, $\delta,$ with the same endpoints as $\gamma'$ stays  less than a distance $2$ from the union of $\tilde{M}_1$ and $\tilde{M}_2.$ Since $\gamma'$ is length-minimizing in this set, it follows that $\gamma'=\delta$ and thus $\gamma'$ does not have any corners. This completes the sketch.

\medskip
\noindent{\bf Claim 1.} The inclusion $incl:M\hookrightarrow Y$  induces a $\pi_1$-isomorphism.

\medskip
  \noindent{\bf Proof of claim 1.}   First we show $incl_*$ is surjective.  Suppose that $\gamma$ is a loop in $Y$ based at a point in $M.$ Then $\gamma=\gamma_1\cdot\gamma_2\cdots\gamma_n$ where each $\gamma_i$ is a path contained  either in $Y_1$ or in $Y_2.$ The endpoints of $\gamma_i$ are in $Y_1\cap Y_2.$ Using condition (6) we can homotop $\gamma$ so that for all $i$ the endpoints of $\gamma_i$ are contained in $M_1\cap M_2.$ Suppose $\gamma_i$ is contained in $Y_j.$ Since  $Y_j$ is a thickening of $M_j$ we may homotop $\gamma_i$ into $M_j$ keeping the endpoints fixed. Thus we may homotop $\gamma$ keeping endpoints fixed into $M.$ Hence $incl_*$ is surjective.
  
Now we show $incl_*$ is injective. By (3) $Y_i$ is a thickening of $M_i$ thus $M_i\hookrightarrow Y_i$ induces a $\pi_1$-isomorphism.  Suppose that $\gamma$ is an essential loop in $M$ which is contractible in $Y.$ By (\ref{convexinjective}) $\gamma$ is not contained in $M_i.$ Thus $\gamma=\gamma_1\cdot\gamma_2\cdots\gamma_n$ where each $\gamma_i$ is a path contained in either $M_1$ or $M_2.$ Using convexity of $M_1$ and $M_2$ we may assume each $\gamma_i$ is a geodesic and has both endpoints in $M_1\cap M_2.$
 We may suppose that $\gamma$ is chosen so that $n$ is minimal. This implies that no $\gamma_i$ is contained in $M_1\cap M_2.$ 

Without loss of generality suppose that $\gamma_i$ is contained in $M_1.$ We claim that $\gamma_i$ is not contained in $Y_2.$ Otherwise $\gamma_i$ is a geodesic in the convex manifold $Y_2$ with both endpoints in the convex submanifold $M_2.$ Consider the universal cover $p:\tilde{Y_2}\rightarrow Y_2.$ Then $\tilde{M}_2=p^{-1}M_2$ is the universal cover of $M_2$ and in particular is connected. A lift, $\tilde{\gamma}_i,$ of $\gamma_i$ is a geodesic in $\tilde{Y}_2$ with both endpoints in $\tilde{M}_2.$ Since $\tilde{M}_2$ is convex it follows that $\tilde{\gamma}_i$ is contained in $\tilde{M}_2$ and thus $\gamma_i$ is contained in $M_2.$ This  contradicts that $\gamma_i$ is not contained in $M_1\cap M_2.$ It follows, that for each $i,$ that $\gamma_i$ is contained in exactly one of $Y_1$ and $Y_2.$

 By (\ref{convexinjective}) $Y_i\hookrightarrow Y$ is $\pi_1$-injective. Thus the universal cover $\tilde{Y}$ of $Y$ is a union of copies of the universal covers of $Y_1$ and $Y_2.$ By (\ref{convexintersect}) the components of $Y_1\cap Y_2$ are convex, thus $\pi_1$-injective into $\pi_1Y$ by (\ref{convexinjective}).  It follows that $\pi_1Y$ is a graph of groups.  By Serre, see \cite{Serre}, the copies of the covers of $Y_1$ and $Y_2$ fit together in a tree-like way to give $\tilde{Y}.$ 
It follows that a lift of the path $\gamma$ to $\tilde{Y}$ does not start and end at the same point.   Hence $\gamma$ is not contractible in $Y.$ This proves $incl_*$ is injective.\qed

\medskip
For $i\in\{1,2\}$ after replacing $M_i$ by its metric completion, we may assume that $M_i$ is a complete metric space. 
The no bumping condition has the following consequence. Suppose that $\alpha$ is an arc in $Y$ which has both endpoints in $M_2\cap\partial M_1.$ If $length(\alpha)\le 2\kappa$ then $\alpha$ is homotopic rel endpoints into $M_2.$ A similar result holds with $M_1$ and $M_2$ interchanged.

\medskip Let $\pi:\tilde{Y}\rightarrow Y$ denote the universal cover and
$dev : \tilde{Y}\rightarrow{\mathbb H}^n$ the developing map. Claim (1) implies that $\tilde{M}=\pi^{-1}M$ is the universal cover of $M.$  As observed above, $\tilde{Y}$ is the union of covering translates of the universal covers, $\tilde{Y}_1$ and $\tilde{Y}_2$ of $Y_1$ and $Y_2.$
 Since the inclusion induces an isomorphism between $\pi_1M_i$ and $\pi_1Y_i$ it follows that each component of $\pi^{-1}M_i$ is a copy of the universal cover, $\tilde{M}_i,$ of $M_i.$ Furthermore every covering translate of $\tilde{Y}_i$ contains exactly one covering translate of $\tilde{M}_i.$   The following claim implies that  the developing map 
$dev:\tilde{M}\rightarrow{\mathbb H}^n$ is injective. The theorem then follows from proposition (\ref{injectivedeveloping}).

 \medskip
 \noindent{\bf Claim 2.} Suppose that $P_0,P_1$ are two points in $\tilde{M}.$ Then there is a hyperbolic geodesic, $\gamma,$ in the interior of $\tilde{Y}$ connecting $P_0$ and $P_1.$ 
 
 \medskip 
  \noindent{\bf Proof of claim 2.}  By condition (2) $M_i$ is convex, so there is a $2$-thickening $M_i^+=T_2(M_i).$ Thus each component, $\tilde{M}_i,$  of $\pi^{-1}(M_i)$ has  a 2-thickening  $\tilde{M}_i^+,$ which, by condition (4), is contained in $\tilde{Y}.$  The covering translates of $\tilde{M}_1$ are pairwise disjoint (and similarly for $\tilde{M}_2$) however the covering translates of $\tilde{M}_1^+$ need not always be disjoint.   For example this will  happen if $M_1$ contains a rank-1 cusp which is contained in a rank-2 cusp of $M.$  Thus the natural isometric immersion of $M_i^+$ into $Y$ is not always injective;  the thickening may bump into itself inside $M_1\cap M_2.$
However the no bumping condition implies that every point of intersection of two different translates of $\tilde{M}_1^+$ is contained in some $\tilde{M}_2$ (and similarly with the roles of $M_1$ and $M_2$ reversed).

 Although $M$ might not be  convex we define $$\tilde{M}^+ = \{\ x\in\tilde{Y}\ :\ d(x,\tilde{M}) \le 2\ \}.$$  This is the union of covering translates of $\tilde{M}_1^+$ and $\tilde{M}_2^+.$     We are assuming that $M_1$ and $M_2$ are complete metric spaces thus $\tilde{M}_1^+,\tilde{M}_2^+$ and $\tilde{M}^+$  are all metrically complete.  
 
 \medskip\noindent{\bf Claim 3.} Suppose $\tilde{M}_i$ is any component of $\pi^{-1}(M_i).$ If $\tilde{M}_1^+\cap\tilde{M}_2^+\ne\phi$ then $\tilde{M}_1\cap\tilde{M}_2\ne\phi.$

\medskip
 \noindent{\bf Proof of claim 3.}  Choose a point $x$ in $\tilde{M}_1^+\cap\tilde{M}_2^+.$ Let $\tilde{Y}_i$ be the component of $\pi^{-1}(Y_i)$ which contains $\tilde{M}_i.$ Since $\tilde{M}_i$ is complete there is a point $a_i\in\tilde{M}_i$ which minimizes $d(a_i,x).$ It follows that $\pi(a_i)\in M_i\setminus int(M_i).$  By definition of $\tilde{M}_i^+$ we have $d(a_i,x) \le 2,$ hence $d(a_1,a_2)\le 4.$
 
 The first case is that $\pi(a_1) \notin int(M_2)$ hence $\pi(a_1)\in M_1\setminus int(M_1\cap M_2).$  Condition (5) implies $N_K(a_1)\subset\tilde{Y}_1.$ Since $\kappa\ge 4$ it follows that $a_2\in\tilde{Y}_1$ thus $a_2\in\tilde{Y}_1\cap\tilde{Y}_2\ne\phi.$ Let $C$ be the component of  $\tilde{Y}_1\cap\tilde{Y}_2$ which contains   $a_2.$ By condition (6) there is a point of $M_1\cap M_2$ in $\pi(C).$ Thus $\pi^{-1}(M_1\cap M_2)$ contains a point, $y,$ in $C.$ However $\pi^{-1}(M_i)\cap\tilde{Y}_i=\tilde{M}_i$  thus $y\in\tilde{M}_1\cap\tilde{M}_2\ne \phi.$

The remaining case is that  $\pi(a_1) \in int(M_2).$ Since $\pi(a_1)\notin int(M_1)$ it follows that $\pi(a_1)\in M_2\setminus int(M_1\cap M_2).$ Let $\tilde{Y}_2'$ be the component of $\pi^{-1}(Y_2)$ which contains $a_1.$ Condition (5) implies $N_K(a_1)\subset\tilde{Y}_2'.$ Since $\kappa\ge 4$ it follows that $a_2\in\tilde{Y}_2'$ thus $\tilde{Y}_2'=\tilde{Y}_2,$ and $a_1\in\tilde{Y}_1\cap\tilde{Y}_2\ne\phi.$ Let $C$ be the component of  $\tilde{Y}_1\cap\tilde{Y}_2$ which contains   $a_1.$ By condition (6) there is a point of $M_1\cap M_2$ in $\pi(C).$ The rest of the argument is the same as the first case. This proves claim 3. \qed
 
\medskip
Since $\tilde{M}^+$ is a complete metric space  there is a length minimizing path
$\gamma:[0,1]\rightarrow\tilde{M}^+$ between the points $P_0$ and $P_1.$

\medskip\noindent{\bf Claim 4.} $\gamma$ is a geodesic except, possibly, at points in $\partial\tilde{M}_1^+\cap\partial\tilde{M}_2^+.$

\medskip
 \noindent{\bf Proof of claim 4.} Consider a point $p$ in the interior of $\gamma.$ If $p$ is in the interior of $\tilde{M}^+$ then, since $\gamma$ is length minimizing, $\gamma$ is a geodesic in a neighborhood of $p.$ Now suppose that $p$ is not in any translate of $\tilde{M}_2^+.$ Thus $p$ is in some $\tilde{M}_1^+.$ Every point of intersection between distinct copies of $\tilde{M}_1^+$ is contained in some $\tilde{M}_2.$ Thus there is an open arc, $\alpha,$ in $\gamma$ that contains $p$ and $\alpha$ is contained in a  unique copy of $\tilde{M}_1^+.$ Since $\tilde{M}_1^+$ is convex, and $\alpha$ is length minimizing, it follows that $\alpha$ is a geodesic. A similar conclusion holds if $p$ is not in any translate of $\tilde{M}_1^+.$ Thus if $\gamma$ is not a geodesic near $p$ then $p$ is in the boundary of $\tilde{M}^+$ and also contained in copies of $\tilde{M}_1^+$ and $\tilde{M}_2^+.$ Hence it is in the boundaries of these copies. This proves claim 4.\qed

\medskip
We will call a point $p$ on $\gamma$ a {\em corner} of $\gamma$ if $\gamma$ is not a hyperbolic geodesic at $p.$ The following claim proves  claim 2 and thus the theorem.

\medskip
\noindent{\bf Claim 5.}   $\gamma$ has no corners. 

\medskip
 \noindent{\bf Proof of claim 5.} If $p$ is a corner of $\gamma$ then (by choosing notation for the covering translates) we may assume  $p\in\partial\tilde{M}_1^+\cap\partial\tilde{M}_2^+.$ It follows that $\tilde{M}_1\cap\tilde{M}_2\ne\phi.$ This is because $\pi(p)\notin M_1$ but $\pi(p)$ is a distance of $2$ from some point  $x\in M_2$ so by condition (5) $\pi(p)\in Y_2.$ Similarly $\pi(p)\in Y_1.$

Let $\delta=\overline{pw}\subset
\tilde{M}_2^+$ denote the maximal subarc of $\gamma$ in $\tilde{M_2}^+$ which
contains $p.$

\centerline{\epsfysize=60mm
\epsfbox{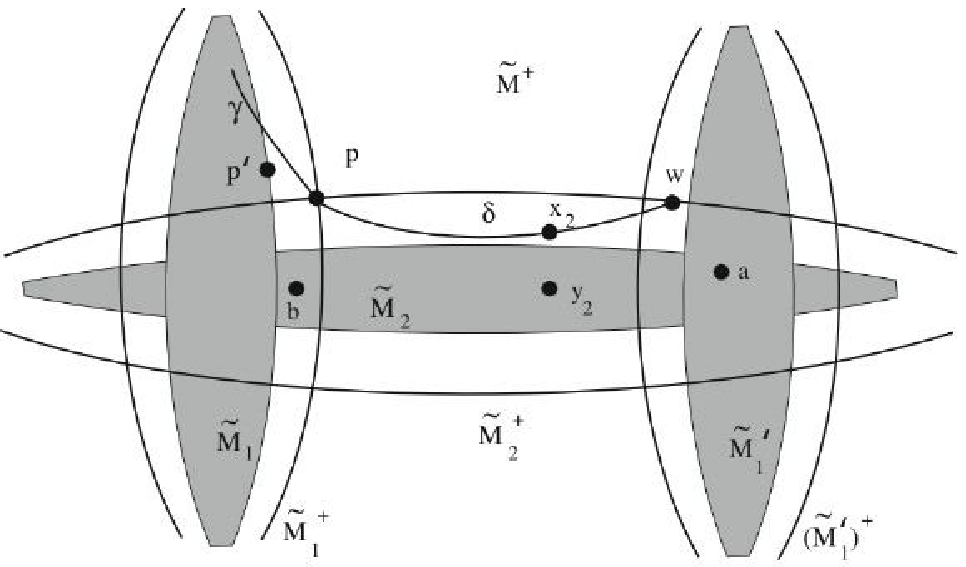}}   

\noindent{\bf Claim 6.} Either $w\in\tilde{M}_2$ or $length(\delta)\ge \kappa-2.$

\medskip \noindent{\bf Proof of claim 6.} First we consider the case that $w$ is an endpoint of $\gamma$ and in addition $\pi(w)\notin M_1.$ Then $\pi(w)\in M_2$ thus $w$ is in some translate, $\tilde{M}_2',$ of $\tilde{M}_2.$ By definition of $\delta$ we have that $w$ is in $\tilde{M}_2^+.$ Thus $\tilde{M}_2'$ intersects $\tilde{M}_2^+$ in the point $w.$ Since $w$ is not contained in any translate of $\tilde{M}_1,$ the no bumping condition implies that $\tilde{M}_2'=\tilde{M}_2.$ Thus $w$ is in $\tilde{M}_2$ and we are done.

Otherwise  $w$ is in some translate, $(\tilde{M}_1')^+,$ of $\tilde{M}_1^+.$  This is because either $w$ is an endpoint of $\gamma$ and $\pi(w)\in M_1$ or else, by maximality of $\delta,$ the point $w$ is in the boundary of $\tilde{M}_2^+$ and thus in some $(\tilde{M}_1')^+$ (since distinct translates of $\tilde{M}_2^+$ can only intersect in a translate of $\tilde{M}_1.$)
Let $\tilde{Y}_1,\tilde{Y}_1'$ be the components of $\pi^{-1}(Y_1)$ which contain $\tilde{M}_1$ and $\tilde{M}_1'$ respectively.

We first consider the case that $\tilde{M}_1\ne\tilde{M}_1'.$  Let $p'$ be the point in $\tilde{M}_1$ closest to $p$ thus $d(p,p')=2.$ If $p'$ was in the interior of $\tilde{M}_2$ then $d(p,\tilde{M}_2)<d(p,p')=2$ which contradicts that $d(p,\tilde{M}_2)=2.$ Hence $p'\in closure(\tilde{M}_1)\setminus int(\tilde{M}_1\cap \tilde{M}_2).$  Since $\tilde{Y}_1$ and $\tilde{Y}_1'$ are not equal they are disjoint. The geodesic $\delta$ has  one endpoint  $p\in \tilde{Y}_1$ and the other endpoint $w\in\tilde{Y}_1'$ thus $w$ is not in $\tilde{Y}_1.$ From  condition (5) (no bumping)  it follows that $d(p',w)\ge \kappa.$ Hence $length(\delta)=d(p,w)\ge \kappa-2.$

The remaining case is that $\tilde{M}_1 = \tilde{M}_1'$ in which case  $\tilde{Y}_1=\tilde{Y}_1'.$ If $\delta$ is contained in $\tilde{Y}_1$ then $\delta$ is a geodesic in the convex set $\tilde{Y}_1$ with both endpoints in the convex subset $\tilde{M}_1^+.$ But this implies that $\delta$ is contained in $\tilde{M}_1^+.$ This in turn means that there is an open interval in $\gamma$ which contains $p$ and is contained in $\tilde{M}_1^+$ and is therefore a geodesic.  This contradicts that $p$ is a corner. Hence $\delta$ contains a point outside $\tilde{Y}_1$ and then as before we obtain $length(\delta)\ge \kappa-2.$ This proves  claim 6.\qed

\medskip
\noindent{\bf Proof of claim 5, resumed.} The geodesic segment $\delta$ has both endpoints, $w,p$ within a distance $2$ of $\tilde{M}_2$ so we may choose points $a,b\in\tilde{M}_2$ with $d(a,w)\le 2$ and $d(b,p)\le 2.$ Since $\tilde{M}_2$ is convex there is a geodesic, $\overline{ab},$ in $\tilde{M}_2.$  If  $length(\delta)\ge \kappa-2 = 4$ it follows from lemma (\ref{quadlemma}) that there is a point $x_2\in\delta$ and a point $y_2\in\overline{ab}$ with $d(x_2,y_2)\le 1$ hence $d(x_2,\tilde{M}_2)\le d(x_2,y_2)\le 1.$ Otherwise $w\in\tilde{M}_2$ and we choose $x_2=y_2=w$ and then $d(x_2,\tilde{M}_2)= d(x_2,y_2)=0\le 1.$

The same argument shows that if $\delta'$ is the maximal segment of $\gamma$ in $\tilde{M}_1^+$ which contains $p$ then there are points $x_1\in\delta'$ and $y_1\in\tilde{M}_1$ with $d(x_1,\tilde{M}_1)\le d(x_1,y_1)\le1.$ Let $\tilde{Y}_i$ be the component of $\pi^{-1}Y_i$ which contains $\tilde{M}_i.$  Using the convexity of $\tilde{Y}_1$ and $\tilde{Y_2}$ it is easy to see that the developing map restricted to $\tilde{Y}_1\cup\tilde{Y}_2$ is an embedding. Thus we may regard  $\tilde{Y}_1\cup\tilde{Y}_2$ as a subset of ${\mathbb H}^n.$ By convexity of the distance function (proposition \ref{distanceconvex}) it follows that every point on the geodesic $\overline{x_1x_2}$ is less than a distance of $1$ from some point on the geodesic $\overline{y_1y_2}.$ Since $p\in(\tilde{M}_1)^+\cap (\tilde{M}_2)^+$ it follows from claim 3 that $\tilde{M}_1\cap\tilde{M}_2\ne\phi.$ Now $\tilde{M}_1,\tilde{M}_2$ are both convex and have non-empty intersection. Since $y_1\in\tilde{M}_1$ and $y_2\in\tilde{M}_2$ it follows from lemma (\ref{convexsetlemma}) that every point on $\overline{y_1y_2}$ is within a distance $\Delta$ of $\tilde{M}_1\cup\tilde{M}_2.$ Thus every point on $\overline{x_1x_2}$ is less than a distance $1+\Delta$ of $\tilde{M}_1\cup\tilde{M}_2.$

\centerline{\epsfysize=60mm
\epsfbox{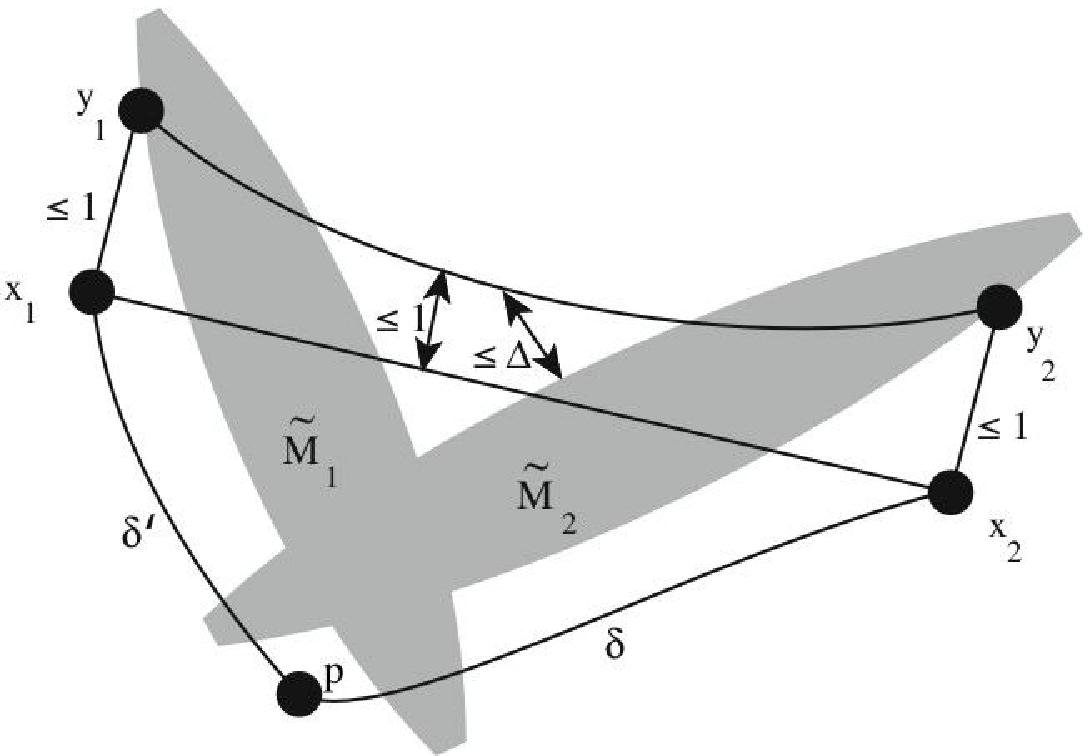}}   

\medskip
The segment, $\delta'\cup\delta,$ of $\gamma$ between $x_1$ and $x_2$ is length minimizing among all paths with the same endpoints in $\tilde{M}_1^+\cup\tilde{M}_2^+.$ Since $1+\Delta < 2$ the hyperbolic geodesic $\overline{x_1x_2}$ is contained in $\tilde{M}_1^+\cup\tilde{M}_2^+$ and is the unique length minimizing path in this set with these endpoints.  Thus  $\delta'\cup\delta=\overline{x_1x_2}$ but this contradicts that $p$ is a corner.  This proves the  claim 5.\qed

\medskip
\noindent{\bf Remark.} Condition (5) (no bumping),  asserts, roughly, that outside $M_1\cap M_2$ that $M_i$ can be $\kappa$-thickened without the thickening bumping into itself. It is this condition which ensures that the  copies of the universal covers of the $M_i$ fit together in hyperbolic space in a treelike way to create $dev(\tilde{M}).$

\begin{lemma}\label{quadlemma} Suppose that $Q$ is a (not necessarily planar) quadrilateral in hyperbolic space with corners $a,b,c,d$ and geodesic sides $\overline{ab},\overline{bc},\overline{cd},\overline{da}.$ Suppose that $|\overline{ad}|\le 2,\ |\overline{bc}|\le 2$  and $|\overline{ab}|\ge4.$ Then there are points  $w\in\overline{ab}$ and $z\in\overline{cd}$   such that $d(w,z)\le 1.$\end{lemma} 
\centerline{\epsfysize=30mm
\epsfbox{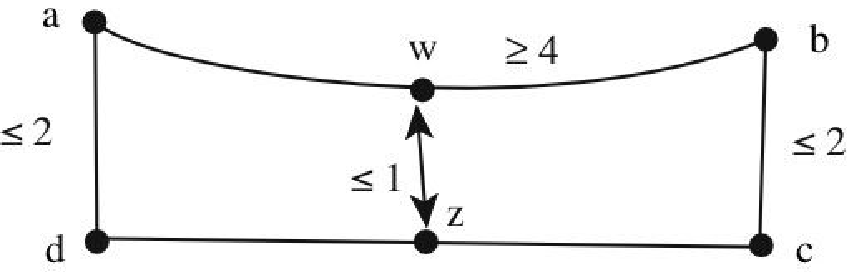}}   
\demo The worst case is the symmetric one in the plane.  A calculation then gives the result.
\qed

\begin{proposition}[distance function is convex \cite{Th1},  p91,  2.5.8]\label{distanceconvex} The distance function $d(x,y),$ considered as a map $d:{\mathbb H}^n\times{\mathbb H}^n\rightarrow {\mathbb R},$ is convex. The composition $d\circ\gamma$ is strictly convex for any geodesic $\gamma$ in ${\mathbb H}^n\times{\mathbb H}^n$ whose projections to the two factors are distinct.
\end{proposition}

\begin{lemma}[convex unions]\label{convexsetlemma} Suppose $A$ and $B$ are convex subsets of hyperbolic space
${\mathbb H}^n$ which have non-empty intersection. Then the convex hull
$CH(A\cup B)$ is contained in the $\Delta$-neighborhood of $A\cup B.$ Here
$\Delta\ =\ log((3+\sqrt{5})/2)\ <1$ is the thin-triangles constant of
${\mathbb H}^2.$
\end{lemma}
\demo Let $X$ be the union of all the geodesic segments, $[a,b]$ with endpoints $a\in A$ and $b\in B.$  We claim  that $X$ is the convex hull of $A\cup B.$ Clearly $X$ is contained in this convex hull. It suffices to show that $X$ is convex. 

Suppose that $p_1\in[a_1,b_1]$ and $p_2\in[a_2,b_2]$ are two points in $X.$ Then we need to show that every point, $q,$ on the geodesic segment from $p_1$ to $p_2$ is also in $X.$
Since $q$ is in the convex hull of the four points $a_1,a_2,b_1,b_2$ there
 is a geodesic $[a,b]$  with endpoints  $a\in [a_1,a_2]\subset A$ and $b\in [b_1,b_2]\subset B$ 
which contains $q.$ Since $[a,b]$ is in $X$ it follows $q\in X$ thus $X$ is convex.

Given $q\in X$ there is $a\in A$ and $b\in B$ such that $q\in [a,b].$   Choose a point $p\in A\cap B$ and consider the
geodesic triangle with sides $[a,p],[p,b],[a,b].$ By convexity
$[a,p]\subset A$ and $[p,b]\subset B.$ The point $q$ on $[a,b]$ is within a
distance $\Delta$ of some point in $[a,p]\cup[p,b]$ and is thus within a
distance $\Delta$ of $A\cup B.$\qed

\section{Cusps in Convex Hyperbolic Manifolds.}
In this section we study the geometry of  cusps in  convex, geometrically-finite, hyperbolic manifolds. 
The main result we need is (\ref{Trel})(e) which states that every thin cusp is contained in a product cusp which is contained in some relative thickening. This is used in the proof of the virtual simple gluing theorem (\ref{virtualsimpleglue}).

If $\Gamma\subset Isom({\mathbb H}^n)$ is  discrete and torsion-free  then the quotient is a geodesically-complete hyperbolic manifold $M={\mathbb H}^n/\Gamma.$ If the limit set of $\Gamma$ contains more than 1 point then, given $\epsilon\ge0,$ we define $C^{\epsilon}$ to be the closed $\epsilon$-neighborhood in ${\mathbb H}^n$ of  the convex hull of the limit set of $\Gamma.$ The {\em $\epsilon$-thickened convex core of $M$} is $Core^{\epsilon}(M)=C^{\epsilon}/\Gamma$ and when $\epsilon=0$  this is called the {\em convex core of $M$}  and we write it as $Core(M).$ The convex core is a convex hyperbolic manifold and a complete metric space. 

The manifold $M$  and the group $\Gamma$ are  {\em geometrically finite} if for some (hence every) $\epsilon>0$ the volume of $Core^{\epsilon}(M)$ is finite. In dimension three this is equivalent to $M$ having a finite sided polyhedral fundamental domain. Bowditch gave several equivalent formulations of geometrical finiteness in \cite{Bow}.

\medskip
{\bf Definitions.} A {\em cusp} is a convex hyperbolic manifold with non-trivial parabolic holonomy. We need to distinguish three kinds of cusp: {\em complete, thin} and {\em product}.  Suppose that $D$ is a closed horoball in ${\mathbb H^n}$ and $\Gamma$ is a non-trivial, discrete, torsion-free subgroup of $Isom({\mathbb H}^n)$ which stabilizes $D.$ Then the quotient $C=D/\Gamma$ is called a {\em complete cusp.} It is a convex hyperbolic manifold with boundary. The induced metric on the boundary of the horoball $D$ is Euclidean, and $\Gamma$ acts by Euclidean isometries on $\partial D,$ so the boundary of $C$ is isometric to a Euclidean $(n-1)$-manifold ${\mathbb E}^{n-1}/\Gamma.$ By a theorem of Bieberbach such a Euclidean manifold is a flat vector bundle over a closed Euclidean manifold.

Suppose that $N={\mathbb H}^n/\pi_1N$ is a geodesically-complete hyperbolic $n$-manifold. A {\em cusp} in $N$ is a submanifold $C^+$ of $N$ which is isometric to a complete cusp. Suppose that $M$ is a convex hyperbolic $n$-manifold and $T_{\infty}(M)={\mathbb H}^n/\pi_1M$ is the corresponding geodesically-complete manifold. Let $C^+$ be a complete cusp in $T_{\infty}(M).$ A {\em cusp} in $M$ is the intersection $C = C^+\cap M.$  We say that $C^+$ is the {\em complete cusp corresponding to $C.$} Clearly $\pi_1C\cong\pi_1C^+.$
   The {\em cusp boundary}  of the cusp $C$ in $M$ is denoted by $\partial_c C$ and equals $(\overline{M\setminus C}) \cap C = C\cap\partial C^+$ and is a submanifold of $\partial C^+.$
   
The horoball $D$ has a codimension-$1$ foliation by horospheres $H_t$ for $t\ge 0$ such that $\partial D=H_0$ and the distance between $H_s$ and $H_t$ is $|s-t|.$ This foliation is preserved by $\Gamma.$ Thus every cusp has a codimension -$1$ foliation whose leaves are called {\em horomanifolds} which are covered by submanifolds of horospheres. The induced metric on a horomanifold is Euclidean.

\medskip
It follows from \cite{Bow} that if $M$ is a convex, geometrically-finite hyperbolic manifold of finite volume and if $C$ is a maximal collection of pairwise disjoint cusps then $\overline{M\setminus C}$ is compact.
 
 \medskip
 \begin{lemma}[thinning cusps]\label{thincusps}  Suppose $M$  is a convex hyperbolic manifold and $C$ is a cusp of $M$ and $\partial_cC\subset \overline{M\setminus C}.$ Then $CH(\overline{M\setminus C}) = (\overline{M\setminus C}) \cup CH(\partial_c C).$
 \end{lemma}
 \demo Since $\partial_cC\subset \overline{M\setminus C}$ it follows that $(\overline{M\setminus C}) \cup CH(\partial_cC) \subset CH(\overline{M\setminus C}).$ Clearly $X \equiv (\overline{M\setminus C}) \cup CH(\partial_cC)$ is closed thus it only remains to show it is also convex.
 
Let $\pi:\tilde{M}\rightarrow M$ be the universal cover. Let $\tilde{X} = \pi^{-1}(X).$ Given two points in $\tilde{X}$ the convexity of $M$ implies there is a geodesic, $\gamma,$ in $\tilde{M}$ connecting them. This geodesic is made up of segments;  each segment is either contained in $\pi^{-1}(\overline{M\setminus C})$ or in a component of $\pi^{-1}(C).$  

Consider a segment, $\delta,$ which is a component  of $\gamma\cap \pi^{-1}(C).$ Each endpoint of $\delta$ is either in $\pi^{-1}(\partial_c C)$ or is an endpoint of $\gamma$ and thus in $\pi^{-1}(CH(\partial_c C)).$  Each component of $\pi^{-1}(C)$ intersects (and contains) a unique component of $\pi^{-1}(CH(\partial_c C)).$ Thus both endpoints of $\delta$ are  contained in the same component, $A,$ of $\pi^{-1}(CH(\partial_c C)).$  Since $A$ is convex $\delta$ is contained in $A.$  Thus $\gamma$ is contained in $\tilde{X}.$ It follows from the definition of convex manifold that $X$ is convex.\qed

 \medskip
{\bf Definition.} A {\em product cusp} is a cusp $C$ such that:\\
(1) $\partial_c C$  is a compact, convex, Euclidean manifold possibly with non-empty boundary, and\\
(2) $C=CH(\partial_c C),$ i.e. $C$ is the convex hull of its cusp boundary.

\begin{proposition}[product cusps are warped products]\label{productcusp} We may isometrically identify ${\mathbb H}^n$ with the upper half space $x_n>0$ of ${\mathbb R}^n$ equipped with the metric $ds/x_n.$ Suppose that $C$ is a product cusp in a hyperbolic $n$-manifold. Then the universal cover of $C$ is isometric to $$\Omega\ =\ \{\ (x_1,\cdots,x_n)\in{\mathbb R}^n\ :\ (x_1,x_2,\cdots,x_{n-1},1)\in S\ and\ x_n \ge 1\ \}$$ 
where $S$ is a convex submanifold of the horosphere $x_n=1$ and $C=\Omega/\pi_1C.$ 
\end{proposition}
\demo Let $D$ be the horoball $x_n\ge1.$ We may isometrically identify the universal cover, $\tilde{C},$ of $C$ with a subset of $D$ such that $S=\tilde{C}\cap\partial D$ is the universal cover of the cusp boundary of $C.$ 

 Since $\tilde{C}$ is convex and limits on the point at infinity it contains $\Omega.$ Since $S$ is contained in $\Omega$ it only remains to show that $\Omega$ is convex. Let $p=(a,t)$ and $p'=(a',t')$ be two distinct points in $\Omega$ with $a,a'\in{\mathbb R}^{n-1}$ and $t,t'\in[1,\infty).$ There is a two-dimensional affine subspace,  $T\subset {\mathbb R}^n,$ which contains $(a,0),(a',0)$ and $(a,1).$ Observe that $T$ is vertical in the sense it contains a line parallel to the $x_n$-axis. Thus $T$ contains $p$ and $p'.$ The  half-plane in $T$ for which $x_n>0$ is a totally geodesic hyperbolic plane in ${\mathbb H}^n$ which contains $p$ and $p'.$ Let $\gamma$ be the unique geodesic arc  in this plane with endpoints $p$ and $p'.$ Let $\pi:{\mathbb R}^n\rightarrow\partial D$ be  given by $\pi(x_1,\cdots,x_{n-1},x_n)=(x_1,\cdots,x_{n-1},1).$ Then $\pi(\gamma)$ is a Euclidean geodesic segment in $\partial D$ with endpoints in $S.$ Since $S$ is convex $\pi(\gamma)$ is contained in $S.$ Since the $x_n$-coordinate of every point on $\gamma$ is at least $1$ it follows that $\gamma$ is in $\Omega.$ Thus $\Omega$ is convex. \qed

 \medskip
The key property of a product cusp that we shall use is that it has a $1$-dimensional foliation by geodesic segments orthogonal to the foliation by horomanifolds. Each geodesic starts on the cusp boundary and goes to infinity. These geodesics are covered by vertical line segments in $\Omega$ which are parallel to the $x_n$-axis.

 \medskip{\bf Definition.} A cusp is {\em thin} if there is a constant $D$ such that for all $t\ge0$ the diameter of the horomanifold $H_t$ in the cusp is less than $D\exp(-t).$ Clearly product cusps are thin.

\begin{lemma}[convex hulls of cusps]\label{cusphull}  Suppose $C$ is a cusp and $\partial_c C$ has bounded diameter. Then $CH(\partial_c C)$ is a thin cusp.\end{lemma}
\demo Let $C^+$ be the complete cusp corresponding to $C.$ Then $\partial_c C$ is a bounded subset of the Euclidean manifold $W=\partial C^+.$ Now $W$ is a vector bundle over a closed Euclidean manifold. Given $r>0$ let $W_r$ be the subset of vectors of length at most $r.$ Then $W_r$ is a compact, convex, submanifold of $W.$ Hence $CH(W_r)$ is a product cusp. Since $\partial_c C$ has bounded diameter for some $r>0$  we have $\partial C\subset W_r.$
Thus $CH(\partial_c C)\subset CH(W_r).$  Product cusps are thin, so $CH(W_r)$ is thin, therefore $CH(\partial_c C)$ is also thin.
\qed

\begin{proposition}[GF cores have thin cusps]\label{flatcusps} Suppose that $N={\mathbb H}^n/\Gamma$ is a  geometrically-finite hyperbolic manifold.  Then $M = Core(N)$ has thin cusps.
\end{proposition} 
\demo If $C$ is a cusp of $M$ then $\partial_c C$ has bounded diameter. This is because otherwise the $\epsilon$-neighborhood of $\partial_c C$ in $N$ has infinite volume. But $N$ is geometrically finite so $C^{\epsilon}(M)$ has finite volume and contains the $\epsilon$-neighborhood of $\partial_c C.$
Since $M=Core(N)$ it follows that $CH(M\setminus C) = M$ and (\ref{thincusps}) implies $V=CH(\partial_c C).$ The result follows from (\ref{cusphull}).
 \qed

 \medskip
Although we won't use this fact, in a geometrically-finite hyperbolic 3-manifold every cusp  contains a (possibly smaller) cusp which is a product cusp. For rank-2 cusps this is obvious. For rank-1 cusps, see \cite{Th2}.

   \medskip
{\bf Definition.} Suppose that $M$ is a convex hyperbolic manifold and that $C_1,\cdots,C_n$ is a maximal collection of pairwise disjoint  cusps in $T_{\infty}(M).$  Given $K\ge0$  the {\em $K$-thickening of $M$ relative to $C=\bigcup_i C_i$} is
 $$T^{rel}_K(M;C)\ =\ CH(T_K(M)\setminus int(C)).$$
Since $M$ is convex, $T_K(M)$ is also convex and, as we shall see below, the difference between it and $T^{rel}_K(M;C)$ is to replace the cusps in $T_K(M)$ by thin cusps. Although relative thickening depends on a particular choice of maximal cusps, usually the choice is unimportant. We will therefore use $T^{rel}_K(M)$ to denote the result of some relative thickening.

\centerline{\epsfysize=60mm
\epsfbox{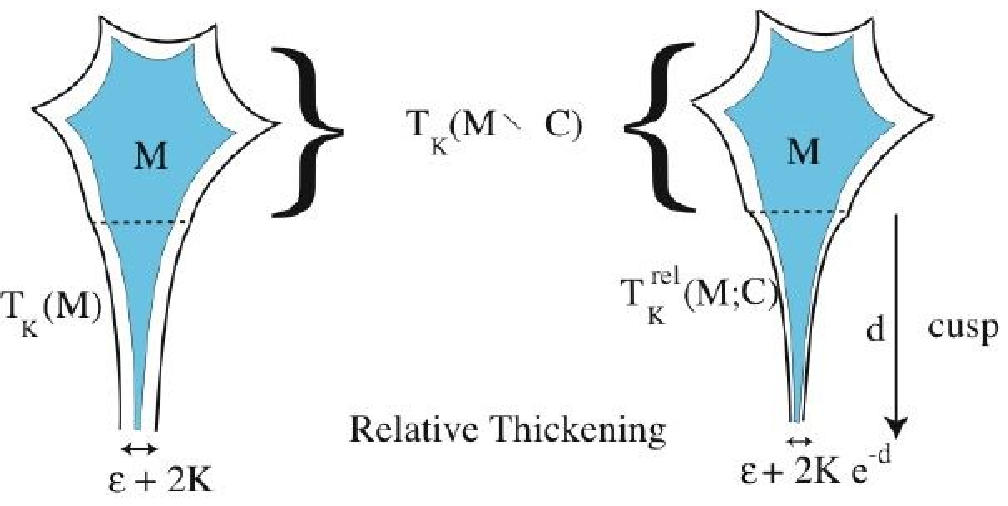}}   

 \medskip
 {\bf Example.} Suppose that $F$ is a complete hyperbolic punctured torus. Isometrically embed ${\mathbb H}^2$ into ${\mathbb H}^3$ then the holonomy of $F$ gives a Kleinian group $\Gamma\cong\pi_1F.$ We can regard  $F$ as a degenerate hyperbolic 3-manifold of zero thickness. The quotient by $\Gamma$ of the $K$-neighborhood of ${\mathbb H}^2$ in ${\mathbb H}^3$ is the $K$-thickening $M=T_K(F).$ It is a convex hyperbolic 3-manifold with a rank-$1$ cusp. The thickness of the cusp everywhere is $2K.$ This means that every point in $F$ is contained a geodesic segment in $M$ of length $2K$ which is orthogonal to $F.$ In particular this is an example of a convex 3-manifold of finite volume which has a rank-1 cusp that is not a product cusp.  The relative $K$-thickening of $F$ is the subset of $T_K(F)$ obtained by replacing the rank-1 cusps of $M$ by product cusps whose thickness decreases exponentially. 
 
\medskip
 \begin{proposition}[relative thickenings contain product cusps]\label{Trel}\  \\ Suppose that $M$ is a convex hyperbolic manifold, and that $C=\bigcup_i C_i$ is a maximal set of pairwise disjoint cusps in $T_{\infty}(M),$ and assume that $M\not\subset C.$ Suppose that each cusp boundary $\partial_c (M\cap C_i)$ in $M$ has bounded diameter.  Then:\\
 (a) $T^{rel}_K(M;C)$  is a convex hyperbolic manifold which contains an
 isometric copy of $M$ and is unique up to isometry fixing $M.$\\
 (b) If $M=Core(M)$ is geometrically finite then  $T^{rel}_K(M;C)$ has finite volume.\\
 (c) If $x\in \overline{M\setminus C}$ and $y\in\partial T_K^{rel}(M;C)$ then $d(x,y)\ge K.$\\
 (d)  $T^{rel}_K(M;C)$ has thin cusps.\\
 (e) For $K$ sufficiently large, each cusp  $M\cap C_i$ is contained in a product cusp which is a subset of $T^{rel}_K(M;C).$
 \end{proposition}
 
 \demo Part (a) is clear. For (b) observe that $T^{rel}_K(M;C)$ has finite volume because it is a subset of the convex manifold $T_K(M),$ and the latter has finite volume because $M$ is geometrically finite. 
For (d) observe that   by (\ref{thincusps}) $$T^{rel}_K(M;C) = (T_K(M)\setminus C) \cup \bigcup_i CH(M\cap\partial_c C_i).$$
Conclusion (d) now follows from (\ref{cusphull}). 

For (c) consider the metric ball, $B,$ of radius $K$ in $T_{\infty}(M)$ centered on $x\in \overline{M\setminus C}.$ Then $B\subset T_K(M).$ Hence $B\setminus T^{rel}_K(M)\subset C.$ Suppose $A$ is a component of $B\cap C.$ We will show that $A\subset T_K^{rel}(M;C).$ It then follows that $B\subset T_K^{rel}(M;C)$ which implies (c).

Identify the universal cover of $T_{\infty}(M)$ with ${\mathbb H}^n.$ Let $\tilde{x}$ be a point in ${\mathbb H}^n$ which covers $x$ and let $\tilde{B}$ be the metric ball in ${\mathbb H}^n$ centered at $\tilde{x}$ and radius $K.$ Thus $\tilde{B}$ projects onto $B.$ Let $C'$ be the component of $C$ which contains $A$ and identify ${\mathbb H}^n$ with the upper half-space model so that the horoball $x_n\ge1$ projects onto $C'.$ Then $\tilde{B}\cap \{\ x_n=1\ \}$ is a metric ball in the horosphere $x_n=1.$ It projects onto $A\cap\partial_cC'.$ The point, $p,$ at infinity ($x_n=\infty$) in the upper half space model is a parabolic fixed point for $M$ and so is in the limit set of $\tilde{M}.$ Let $Y$ be the subset of ${\mathbb H}^n$ corresponding to the universal cover of $T_K^{rel}(M;C).$ It follows that $p$ is in the limit set of $Y.$ Now $Y$ is convex and contains contains $\tilde{B}\cap \{\ x_n=1\ \}$ and limits on $p.$ Hence $Y$ contains the solid cylinder, $W,$ of points lying vertically above  $\tilde{B}\cap \{\ x_n=1\ \}.$ Since $x$ is not in the interior of $C'$ it follows that $\tilde{x}$ is not in $\{\ x_n>1\ \}.$ Hence $W$ contains $Z = \tilde{B}\cap  \{\ x_n\ge1\ \},$ and $Z \subset Y.$ The projection of $Z$ to $T_K^{rel}(M;C)$ is $A.$ This proves (c).

For (e) the cusp boundary $F_i = M\cap C_i$ has bounded diameter. Thus $F_i$ contained in a compact convex submanifold $E_i$ of $\partial C_i.$ For $K$ sufficiently large $E_i$ is contained in the closure of $T_K(M)\setminus C.$ Since $F_i\subset E_i$ it follows that $CH(F_i)\subset CH(E_i)\subset CH(T_K(M)\setminus C).$ Now $CH(E_i)$ is a product cusp and this proves (e). \qed

  \section{Induced  gluing}
The main result of this section is the virtual simple gluing theorem (\ref{virtualsimpleglue}).  Suppose that we have two locally-isometric immersions of  hyperbolic $3$-manifolds equipped with basepoints into a hyperbolic $3$-manifold $M$ $$f:(A,a_0)\rightarrow (M,m_0)  \qquad\qquad g:(B,b_0)\rightarrow (M,m_0).$$  We would like to use this information to glue $A$ and $B$ together. For example if both immersions are injective then we might identify $A$ with $f(A)$ and $B$ with $g(B).$ The union $f(A)\cup g(B)\subset M$ may then be regarded as a quotient space of the disjoint union of $A$ and $B$ and it has a hyperbolic metric. 
 
 However we will be interested in situations when the immersions are not injective. Furthermore, even when $A$ and $B$ are submanifolds of $M,$ we want to do make the fewest identifcations subject to the requirements that the basepoints in $A$ and $B$ are identified and that the identification space is a hyperbolic $3$-manifold. Thus if $A\cap B$ is not connected we wish to only identify $A$ and $B$ along the component, $C,$ of $A\cap B$ containing the basepoint $m_0.$ In certain circumstances the fundamental group of the identification space will be a free product of the fundamental groups of $A$ and $B$ amalgamated along a subgroup corresponding to the fundamental group of $C.$

 We give below a very general way of forming an identification space.  Even when  the resulting identification space is a hyperbolic manifold,  it will usually not be convex and might not have a convex thickening.
 
\medskip
\noindent{\bf Definitions.} Suppose that $f:(X,x_0)\rightarrow (Z,z_0)$ and $g:(Y,y_0)\rightarrow (Z,z_0)$ are
continuous maps of pointed spaces. Define a relation $R$ on the disjoint union, $X\coprod Y,$ of $X$ and $Y$
as follows. If $x\in
X$ and $y\in Y$ then $xRy$ iff there are paths $\alpha:I\rightarrow X$ and
$\beta:I\rightarrow Y$ such that
$\alpha(0)=x_0,\ \alpha(1)=x,\ \beta(0)=y_0,\ \beta(1)=y$ and
$f\circ\alpha=g\circ\beta.$ 
We now define a topological space called the {\em induced gluing} denoted
$S(f,g)$ to be the quotient space $X\coprod Y/ \equiv$ obtained by
taking the equivalence relation $\equiv$ which is generated by $R.$ 
It is clear that if two points are identified then they have the same image in $Z.$ Let $\pi_X:X\rightarrow S(f,g)$ and $\pi_Y:Y\rightarrow S(f,g)$ denote the natural
projections. We say that the induced gluing is a {\em simple gluing} if both $\pi_X$ and $\pi_Y$ are injective.

\medskip
{\bf Example (1)} If $f,g$ are both embeddings define $Z_0$ to be the path component of $f(X)\cap
g(Y)$ containing the basepoint $z_0.$
Then $S(f,g)$ is  obtained from $X\coprod Y$ by identifying $f^{-1}(Z_0)$ with
$g^{-1}(Z_0)$ using the homeomorphism
$g^{-1}\circ f:f^{-1}(Z_0)\rightarrow g^{-1}(Z_0).$ It is clear that this is a simple gluing. 
 
\medskip
{\bf Example (2)} Suppose  $f,g:S^1\times[0,2]\rightarrow S^1\times S^1$ are given by $$f(\omega,t)=(\
\omega\ ,\ exp(2\pi i t)\ )\qquad\text{and}\qquad g(\omega,t)=(\ exp(2\pi i t)\ ,\ \omega\ ).$$
Then $S(f,g)$ can be naturally identified with the codomain $S^1\times S^1.$ In this case the gluing is not simple, in fact both projections are surjective.

\medskip
{\bf Example (3)} Suppose now that
$p:S^1\times [0,2]\rightarrow S^1\times [0,2]$ is the 3-fold cyclic cover and $f,g$ are as in example (2). Then $S(f\circ p,g\circ p)$
is homeomorphic to a torus minus an open disc. It is easy to see that this is a simple gluing. Futhermore this is a modification of example (2) where the domains are replaced by certain finite covers. This phenomenon of taking a non-simple gluing and making a simple gluing by replacing the spaces by finite covers is generalized below.
 
 \medskip
{\bf Example (4)} Suppose that $X= S^1$ and $Y\cong Z\cong D^2$ and $f(z)=z^2$ and $g$ is a homeomorphism. Then $\pi_Y$ is a homeomorphism and $\pi_X$ is a covering map onto $\partial D.$ The gluing is not simple. Furthermore there  are no finite covers of the domains which result in a simple gluing as in example (3). 
 
\medskip It is routine to check the following:
 
\begin{lemma}[induced gluing]\label{simplegluelemma} Suppose  $f:(X,x_0)\rightarrow (Z,z_0)$ and
$g:(Y,y_0)\rightarrow (Z,z_0)$ yield an induced gluing $S(f,g).$\\
(a) The induced gluing is simple iff for every $x\in X$ there is at most one $y\in Y$ with $xRy$ and
vice versa.\\ 
(b) If the induced gluing is simple define subspaces $X_0\subset X,\ Y_0\subset Y$ by $x\in X_0$ and $y\in
Y_0$ if $xRy.$ Define $h:X_0\rightarrow Y_0$ by $h(x)=y$ if $xRy.$ Then $S(f,g)$ is the
quotient space obtained from $X\coprod Y$ by identifying $X_0$ with $Y_0$ using $h.$\\
(c) If $X,Y,Z$ are smooth $n$-manifolds with boundary and $f,g$ are immersions, and if the induced gluing is simple, and  if $f|\partial X$ is transverse to $g|\partial Y,$ then $S(f,g)$ is an
$n$-manifold with boundary.\\
 (d) There is a unique continuous induced map $h:S(f,g)\rightarrow Z$ such that $f=h\circ\pi_X$ and $g=h\circ\pi_Y.$\\
 (e) If the induced gluing is simple then $X\cap Y$ is path connected subspace of $S(f,g).$
\end{lemma}
 
 \medskip
 We are concerned with induced gluings of convex hyperbolic manifolds. In this case the gluing  of the manifolds is determined by the intersection of the  images of their developing maps:
 \begin{lemma}[geodesic gluing]\label{convexinducedglue} Suppose $A,B,M$ are convex hyperbolic manifolds of the same dimension and $f:(A,a_0)\rightarrow (M,m_0)$ and $g:(B,b_0)\rightarrow (M,m_0)$ are locally isometric immersions. Let $p_A:\tilde{A}\rightarrow A$ and $p_B:\tilde{B}\rightarrow B$ and $p_M:\tilde{M}\rightarrow M$ be the universal covers. Choose lifts of the base points and maps $\tilde{f}:(\tilde{A},\tilde{a}_0)\rightarrow(\tilde{M},\tilde{m}_0)$ and $\tilde{g}:(\tilde{B},\tilde{b}_0)\rightarrow(\tilde{M},\tilde{m}_0)$ covering $f$ and $g.$  Suppose $a\in A$ and $b\in B,$ then $aRb$ iff there are $\tilde{a}\in\tilde{A}$ and $\tilde{b}\in\tilde{B}$ covering $a$ and $b$ respectively such that $\tilde{f}(\tilde{a})=\tilde{g}(\tilde{b}).$ Furthermore, if $aRb$ then the paths $\alpha,\beta$ used in the definition of the $R$-relation can be chosen to be geodesics.  \end{lemma}
 \centerline{\epsfysize=50mm
\epsfbox{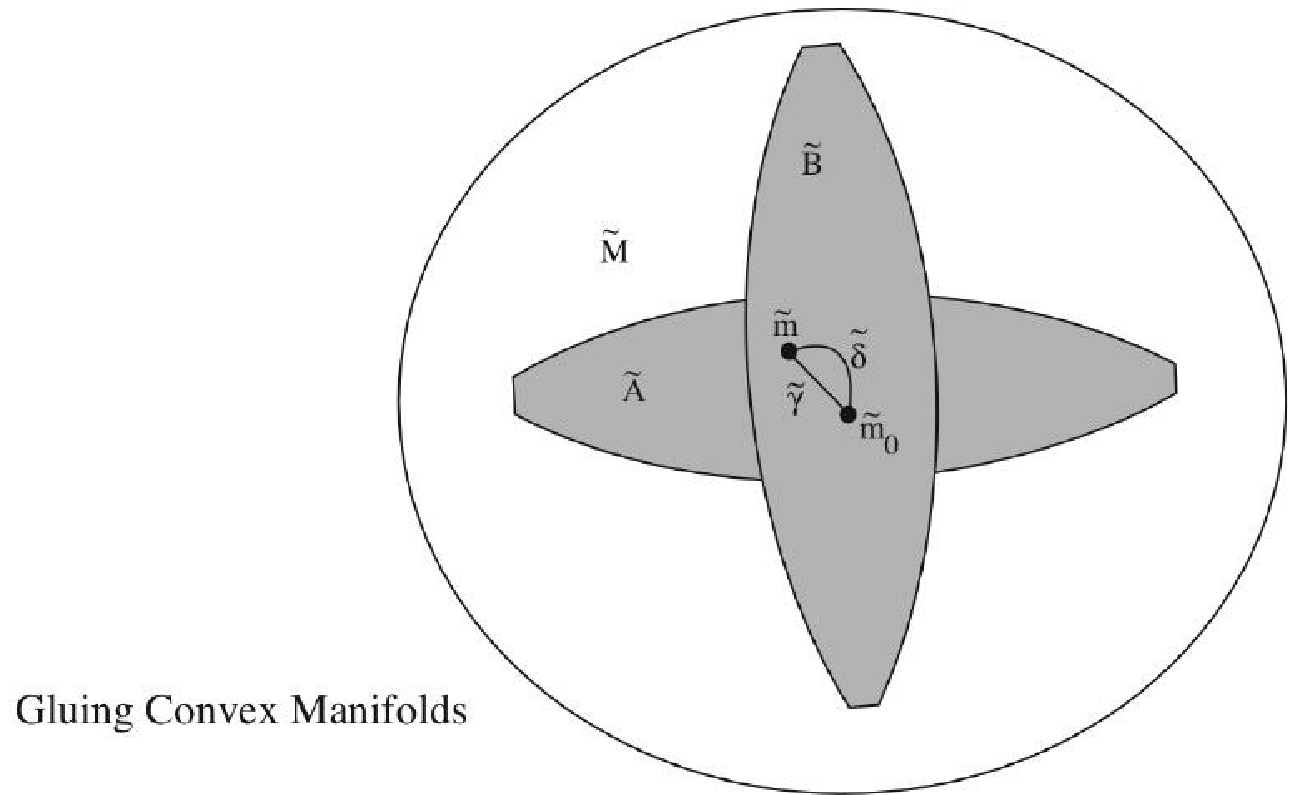}}   
\demo First suppose that $aRb.$ Then there are paths $\alpha$ in $A$ from $a_0$ to $a$ and $\beta$ in $B$ from $b_0$ to $b$ with $f\circ\alpha=g\circ\beta.$ Let $\tilde{\alpha},\tilde{\beta}$ be the lifts that start at the respective basepoints $\tilde{a}_0,\tilde{b}_0$ of $\tilde{A},\tilde{B}.$ Then $\tilde{\delta}=\tilde{f}\circ\tilde\alpha$ is a lift of $f\circ\alpha$ and $\tilde{g}\circ\tilde{\beta}$ is a lift of $g\circ\beta$ which both start at $\tilde{m}_0$ thus $\tilde{f}\circ\tilde\alpha=\tilde{g}\circ\tilde{\beta}.$  Setting $\tilde{a}=\tilde{\alpha}(1)$ gives $p_A(\tilde{a})=a.$ Similarly $\tilde{b}=\tilde{\beta}(1)$ implies $p_B(\tilde{b})=b.$ Since $\tilde{f}\circ\tilde\alpha=\tilde{g}\circ\tilde{\beta}$ we get $\tilde{f}(\tilde{a})=\tilde{g}(\tilde{b})$  completing the proof in this direction.
     
 For the converse, given $\tilde{a},\tilde{b}$ with $\tilde{f}(\tilde{a})=\tilde{g}(\tilde{b}),$ the point $\tilde{m}=\tilde{f}(\tilde{a})$ is in both $\tilde{f}(\tilde{A})$ and $\tilde{g}(\tilde{B}).$ These are convex subsets of $\tilde{M}$ each of which may be identified with a convex subset of hyperbolic space; thus their intersection is a convex set, $C,$ that contains both $\tilde{m}_0$ and $\tilde{m}.$ By convexity there is a geodesic $\tilde{\gamma}$ in $C$ with endpoints $\tilde{m}$ and $\tilde{m}_0.$ Then $\tilde{\alpha}=\tilde{f}^{-1}\tilde{\gamma}$ and $\tilde{\beta}=\tilde{g}^{-1}\tilde{\gamma}$ are geodesics in $\tilde{A}$ and $\tilde{B}$ which project to geodesics $\alpha=p_A\tilde{\alpha}\subset A$ and $\beta=p_B\tilde{\beta}\subset B.$ Set $a=\alpha(1)$ and $b=\beta(1)$ then since $\tilde{f}\circ\tilde{\alpha}=\tilde{g}\circ\tilde{\beta}$ it follows that $f\circ\alpha=g\circ\beta$ and hence $aRb.$ \qed

\medskip
We now generalize the passage from example (2) to example (3). Example (4) shows that the convexity hypothesis is necessary.

\begin{theorem}[virtual simple gluing theorem]\label{virtualsimpleglue} Suppose that $M={\mathbb H}^n/\pi_1M$ is a geodesically-complete hyperbolic manifold. Suppose that $A$ and $B$ are geometrically-finite, convex, hyperbolic $n$-manifolds with finite
volume and thin cusps. Suppose
$f:(A,a_0)\rightarrow (M,m_0)$ and $g:(B,b_0)\rightarrow (M,m_0)$ are local isometries. 
Then there are finite covers $p_A:(\tilde{A},\tilde{a}_0)\rightarrow (A,a_0)$ and
$p_B:(\tilde{B},\tilde{b}_0)\rightarrow (B,b_0)$ and maps $\tilde{f}=f\circ p_A:(\tilde{A},\tilde{a}_0)\rightarrow
(M,m_0)$ and $\tilde{g}=g\circ p_B:(\tilde{B},\tilde{b}_0)\rightarrow
(M,m_0)$  such that $S(\tilde{f},\tilde{g})$ is a simple gluing of $\tilde{A}$ and $\tilde{B}.$ 
 
 Furthermore, let
$G=f_*\pi_1(A,a_0)\cap g_*\pi_1(B,b_0) <
\pi_1(M,m_0)$ and suppose that 
$G_A=f_*^{-1}G$ is separable in $\pi_1(A,a_0)$ and $G_B=g_*^{-1}(G)$ is separable in $\pi_1(B,b_0).$
Then the covers can be chosen so that
$G_A\subset p_{A*}(\pi_1(\tilde{A},\tilde{a}_0))$ and $G_B\subset p_{B*}(\pi_1(\tilde{B},\tilde{b}_0)).$
\end{theorem}
\demo  First we give a sketch of the proof, the details follow. In claim 1 we show that if $A$ and $B$ are compact there is a constant $L>0$ such that if  $a_1\in A$ is $R$-related  to $b_1\in B$ then there are two geodesics, one in $A$ and the other in $B,$ each of length at most $L,$ which are identified by $f$ and $g,$ and that start at the basepoint and end at $a_1$ and $b_1.$ 

Using separability there are finite covers $\tilde{A},\tilde{B}$  of $A$ and $B$ such that the only loops of length at most $2L$ which lift to these covers correspond to elements of $G.$ If the extra hypothesis of separabilty is not assumed we set $G=1$ and use residual finiteness. In claim 2 we show that the induced gluing of these covers is simple.
 
Indeed, suppose two points $\tilde{a}_1,\tilde{a}_2$ in $\tilde{A}$ are $R$-related to the same point $\tilde{b}$  in $\tilde{B}.$ Then there are geodesics in $\tilde{A},\tilde{B},$  possibly very long, which are identified. Projecting these geodesics into $A$ and $B$ we will find new geodesics in $A$ and $B$  with the same endpoints as the originals, that are identified, and have length at most $L.$ This pair of geodesics form a loop of length at most $2L.$ We now argue these geodesics lift to the covers and have the same endpoints as the originals.  One deduces $\tilde{a}_1=\tilde{a}_2.$ Thus there are no self-identifications.
 
Finally,  in the non-compact case, we first replace $A$ and $B$ by relative thickenings of $A$ and $B$ so that each thin cusp of $A$ or $B$ is contained in a product cusp which in turn is contained in a cusp of the thickenings. Then we truncate along cusps to obtain compact submanifolds of the relative thickenings.  This allows us to find a constant $L$ that applies to points in these submanifolds and the previous argument shows there are no self-identifications in the submanifolds. The fact that thin cusps of $A$ and $B$ are contained in product cusps means any self-identifications within a thin cusp propagate vertically within  the larger product cusp in a product like way to give identifications where the cusp meets the compact submanifold. Thus there are no self-identifications in the thin cusps of $A$ or $B$ either. This completes the sketch.
   
 \medskip
The first step is to replace $A$ and $B$ by thickenings so that the thin cusps of $A$ and $B$ are contained in product cusps of the thickenings. We start by renaming $A$ as $A^*$ and $B$ as $B^*.$ Let $C$ be a maximal collection of pairwise disjoint cusps in $M$. Then $C_A = f^{-1}(C)$ and $C_B = g^{-1}(C)$ is a maximal collection of pairwise disjoint cusps in $A^*$ and $B^*$ respectively. We choose the cusps in $C$ small enough so that $A^*\ne C_A$ and $B^*\ne C_B$ and so the cusps $C_A$ and $C_B$ are thin. By (\ref{Trel})(e) there is $K>0$ such that each  component of $C_A$ and of $C_B$ is contained in a product cusp contained in $A \equiv T^{rel}_K(A^*)$ or $B \equiv T^{rel}_K(B^*)$ respectively. The maps $f,g$ have natural extensions to local isometries $f:A\rightarrow M$ and $g:B\rightarrow M.$ Henceforth we shall use  $f$ and $g$ denote these extensions. 

Define $A^-=A\setminus f^{-1}(C),$ and  $B^-=B\setminus g^{-1}(C).$ Then $A^-,B^-$ are compact.  Furthermore each component of $C_A$ is conatined in a product cusp contained in $A\setminus A^-$ and similarly for $B.$

\medskip{\bf Claim 1.} There is $L>0$ such that if $a_1\in A^-$ and $b_1\in B^-$ and $a_1Rb_1$  then there are geodesics $\alpha\subset A$ and $\beta\subset B$ with $length(\alpha)=length(\beta)\le L$ and $f\circ\alpha=g\circ\beta$ and $\alpha$ has endpoints $a_0,a_1$ and $\beta$ has endpoints $b_0,b_1.$
 
\medskip \noindent{\bf Proof of claim 1.} Define  $$S = \{ (a_1,b_1)\in A\times B: a_1Rb_1\}\qquad and\qquad \ell:S\rightarrow{\mathbb R}$$ by $\ell(a_1,b_1)$ is the length of some shortest geodesic, $\alpha,$ as above. By (\ref{simplegluelemma})(c) $S(f,g)$ is Hausdorff  thus $S$ is a closed subset of $A\times B.$ We first show that $\ell$ is continuous at points of $S$ in the interior of $A\times B.$
 
  If  $(a_1,b_1)\in S$ and  $a_1$ is  in the interior of $A$ and $b_1$ is in the interior of $B$ then $\ell$ is continuous at $(a_1,b_1).$ This is because there are small open balls  $U\subset A$ and $V\subset B$ centered on $a_1$ and $b_1$ with $fU=gV.$ It is clear that each $a\in U$ is $R$-related to a unique point $b\in V.$ Since $A$ and $B$ are convex there are geodesics, $\alpha',\beta'$ close to $\alpha$ and $\beta$ from the basepoints to $a$ and $b$ and $f\alpha'=g\beta'.$ This proves continuity of $\ell$ at interior points. 
  
  By enlarging $A$ and $B$ to slightly larger convex manifolds, $A^+,B^+$ every point in $S$ is in the interior of $A^+\times B^+$ thus the corresponding function $\ell^+$ defined on $S^+$ is continuous at every point in $S.$ The function $\ell^+$ is defined using geodesics in the enlarged manifolds. However if a geodesic in $A^+$ starts and ends in $A$ then, since $A$ is convex, the geodesic is contained in $A.$ It follows that $\ell^+|S=\ell$ and thus $\ell$ is continuous on all of $S.$ 
  
 Restricting the continuous function $\ell$ to the compact set $S\cap(A^-\times B^-)$ we obtain the required bound $L$ on $\ell.$ 
 Observe that the shortest geodesics used above connecting a point in $A^-$ to the basepoint $a_0$ are not necessarily contained in $A^-,$ they may go into the cusp of $A$ some bounded distance. This proves claim 1.\qed
  
\medskip     Let
$S\subset\pi_1(B,b_0)$ be the set of elements represented  by based loops of length at most $2L.$ Since $B$ is convex $S$ is finite.  To prove the theorem with the additional hypothesis: since $A,B,M$ are all convex the maps $f_*,g_*$ are injective by proposition (\ref{convexinjective}), and for notational simplicity we will
identify $G$ with $G_A$ and $G_B.$
Since $G$ is separable in $\pi_1(B,b_0)$ there is a subgroup, $H < \pi_1(B,b_0),$ of finite index which
contains
$G$ and contains no element of  $S\setminus G.$

Without the additional hypothesis  we set $G=G_A=G_B=1.$ Since $A$ and $B$ are convex hyperbolic $n$-manifolds their fundamental groups are isomorphic to subgroups of $O(n,1)$ and are thus linear and therefore residually finite by \cite{Mal}. This means that the trivial group is separable in $\pi_1(A,a_0)$ and $\pi_1(B,b_0).$ Thus there is a subgroup, $H < \pi_1(B,b_0),$ of finite index which
contains no element of  $S\setminus \{1\}.$
 
  Let $p_B:\tilde{B}\rightarrow B$ be the cover corresponding to this
subgroup. Similarly define $p_A:\tilde{A}\rightarrow A.$ Let $\tilde{f}=f\circ p_A$ and $\tilde{g}=g\circ p_B.$ Define $\tilde{A}^-=p_A^{-1}(A^-)$ and  $\tilde{B}^-=p_B^{-1}(B^-).$ Then $A^-,B^-$ are compact.  
 The following claim and lemma (\ref{simplegluelemma})(a) imply the induced gluing, $S(\tilde{f},\tilde{g}),$ of $\tilde{A}^-$ and $\tilde{B}^-$ is simple.
  
   \medskip{\bf Claim 2.} If $\tilde{a}_1,\tilde{a}_2\in \tilde{A}^-$ are both
$R$-related to $\tilde{b}_1\equiv\tilde{b}_2\in \tilde{B}^-$ then $\tilde{a}_1=\tilde{a}_2$ and similarly with the roles of $\tilde{A}^-$ and $\tilde{B}^-$ reversed.

\noindent{\bf Proof of claim 2.}  Observe that $a_i=p_A(\tilde{a}_i)\in A^-$ and $b_i=p_B(\tilde{b}_i)\in B^-.$
 
\centerline{\epsfysize=60mm
\epsfbox{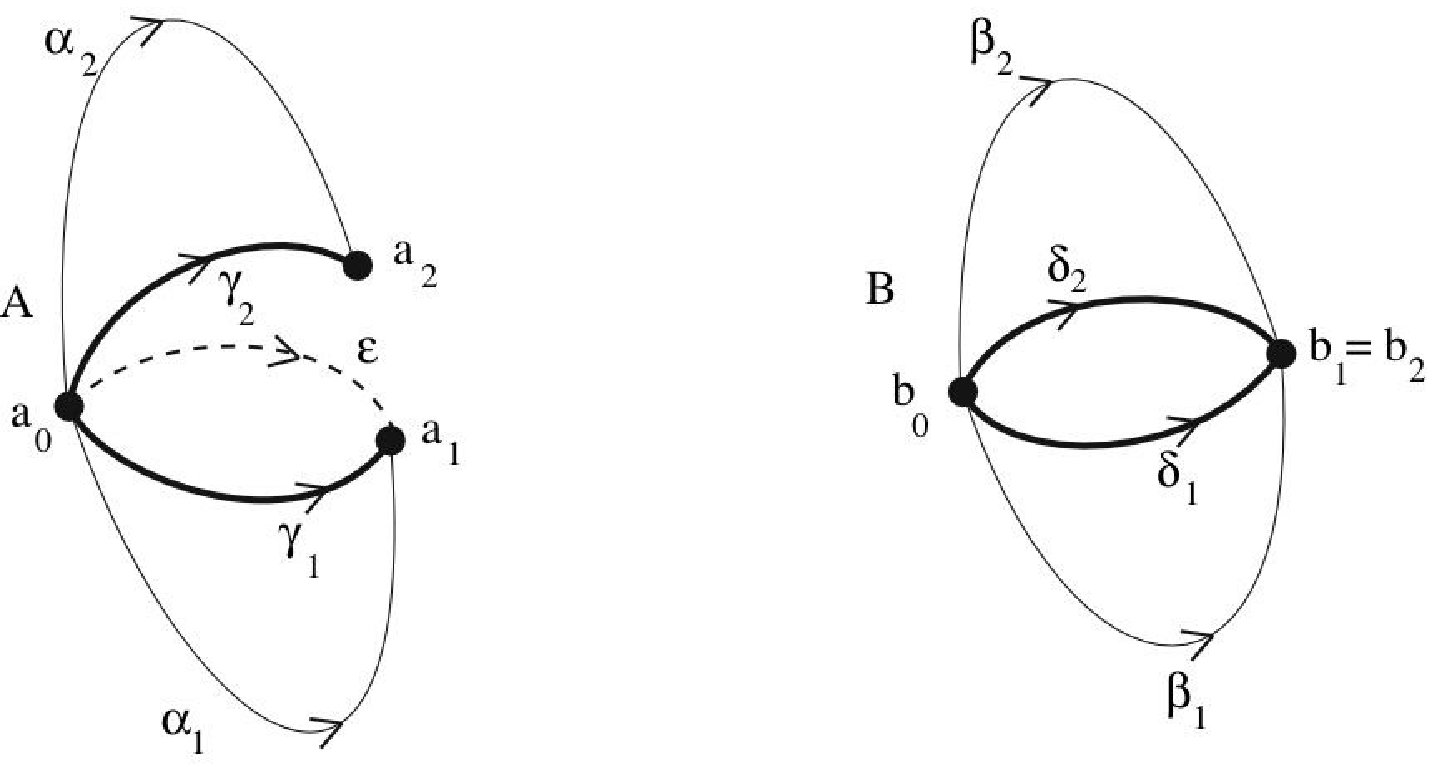}}   
    
For $i\in\{1,2\}$ because $\tilde{a}_iR\tilde{b}_i$ there is a geodesic $\tilde{\alpha}_i\subset\tilde{A}$ starting at $\tilde{a}_0$ and ending
at $\tilde{a}_i$ and a geodesic $\tilde{\beta}_i\subset\tilde{B}$ starting at $\tilde{b}_0$ and
ending at $\tilde{b}_i$ such that $\tilde{f}\circ\tilde{\alpha}_i=\tilde{g}\circ\tilde{\beta}_i.$ Project
these into $A$ and $B$ to obtain geodesics $\alpha_i=p_A\circ\tilde{\alpha}_i$ in $A$ and 
$\beta_i=p_B\circ\tilde{\beta}_i$ in $B.$ Observe that $f\circ\alpha_i=g\circ\beta_i.$ Thus there are geodesics
$\gamma_i\subset A$ and $\delta_i\subset B$  of length at
most $L$ with $f\gamma_i=g\delta_i$ such that the endpoints of $\alpha_i$ and $\gamma_i$ are the same and the endpoints of $\beta_i$ and $\delta_i$ are the same as shown in the diagram.

Consider the loop $\alpha_i\cdot\gamma_i^{-1}$ in $A.$ The image of this loop under $f$ is the same as
the image under $g$ of the loop $\beta_i\cdot\delta_i^{-1}$ hence these loops give elements of $G.$
 It follows that the loops $\alpha_i\cdot\gamma_i^{-1}$
lift to loops in $\tilde{A}$ based at $\tilde{a}_0.$  Similarly  the loops $\beta_i\cdot\delta_i^{-1}$
lift to loops in $\tilde{B}$ based at $\tilde{b}_0.$

\medskip

{\bf Claim 3.} The loop $\delta=\delta_1\cdot\delta_2^{-1}$  lifts to a loop $\tilde{\delta}$  in $\tilde{B}$ based at $\tilde{b}_0.$

\noindent  This is because 
$$\delta = \delta_1\cdot\delta_2^{-1} = (\delta_1\cdot\beta_1^{-1})\cdot(\beta_1\cdot\beta_2^{-1})\cdot(\beta_2\cdot\delta_2^{-1}).$$
The loops $\delta_1\cdot\beta_1^{-1}$ and $\beta_2\cdot\delta_2^{-1}$ lift to loops based at $\tilde{b}_0$ by the previous paragraph. The path $\tilde{\beta}_i$ starts at $\tilde{b}_0$ and ends at $\tilde{b}_i.$ Using the assumption that $\tilde{b}_1=\tilde{b}_2$ we see that the $\tilde{\beta}_1,\tilde{\beta}_2$ have the same endpoints thus $\beta_1\cdot\beta_2^{-1}$ lifts to a loop based at $\tilde{b}_0.$ Thus all three loops in the product lift to loops based at $\tilde{b}_0.$ This proves claim 3.\qed

\medskip
We continue with the proof of claim 2. In what follows $\simeq$ denotes homotopy between maps of an interval keeping endpoints fixed. Since $length(\delta)\le 2L,$ and using the definition of $\tilde{B},$ we see that $[\delta]\in G.$ Since $[\delta]\in G$ there is a loop, $\eta,$ in $A$ based at $a_0$ such that $f\circ\eta\simeq g\circ\delta.$ By sliding one endpoint of $\eta$ along $\gamma_1$ one obtains a path $\epsilon$ in $A$ with endpoints $a_0,a_1$ such that $f\circ(\gamma_1\cdot\epsilon^{-1})\simeq g\circ\delta.$ By changing basepoints it follows that $$f\circ(\epsilon^{-1}\cdot\gamma_1)\simeq g\circ(\delta_2^{-1}\cdot\delta_1).$$ Since $A$ is convex we may homotop $\epsilon$ keeping its endpoints fixed to be a geodesic. 

Combining this with: 
$$\begin{array}{rcl} f\circ(\epsilon^{-1}\cdot\gamma_1) & = & (f\circ\epsilon^{-1})\cdot(f\circ\gamma_1)
\\
g\circ(\delta_2^{-1}\cdot\delta_1) & =  & (g\circ\delta_2^{-1})\cdot(g\circ\delta_1)\\
f\circ\gamma_1 & = & g\circ\delta_1 \end{array}$$ 
 it follows that $f\circ\epsilon^{-1}\simeq g\circ\delta_2^{-1}$ and thus $f\circ\epsilon\simeq g\circ\delta_2.$ Since $f\circ\epsilon$ and $g\circ\delta_2$ are geodesics in the convex hyperbolic manifold $M$ with the same endpoints, and are homotopic rel endpoints, it follows they are equal. But $f\circ\gamma_2$ is also equal to this geodesic. Thus $\gamma_2$ and $\epsilon$ are two geodesics in $A$ which start at the same point and map under the local isometry $f$ to the same geodesic. It follows that $\epsilon=\gamma_2.$

We now have a loop $\gamma_1\cdot\gamma_2^{-1}\subset A$ such that   $f\circ(\gamma_1\cdot\gamma_2^{-1})\simeq g\circ\delta.$ Since the latter is in $G$ it follows that $\gamma_1\cdot\gamma_2^{-1}$ lifts to a loop based at $\tilde{a}_0$ and therefore $\tilde{a}_1=\tilde{a}_2.$  This completes the proof claim 2.\qed

 \medskip We now finish the proof of the theorem. We have shown that the induced gluing of $\tilde{A}^-$ and $\tilde{B}^-$ is simple. Suppose the induced gluing of the corresponding covers $\tilde{A}'$ and $\tilde{B}'$ is not simple. Since $\tilde{A}'\subset \tilde{A}$ it follows that if two points, $a_1,a_2\in \tilde{A}'$ are identified by the gluing then these points are in the cusps $\tilde{C}_A$ of $\tilde{A}'$ that cover $C_A.$ Clearly a finite cover of a product cusp is also a product cusp.  There are product cusps contained in $\tilde{A}$ which contain $\tilde{C}_A.$ A product cusp has a $1$-dimensional foliation by rays starting on the cusp boundary. Hence two of these rays are identified by the induced gluing. This implies two point on the cusp boundaries are also identified. The cusp boundaries are in $\tilde{A}^-$ thus there are two points  in $\tilde{A}^-$ which are identified and this is a contradiction. Similar remarks apply to $B'.$  Hence the induced gluing of $\tilde{A}'$ and $\tilde{B}'$ is simple. This completes the theorem. \qed

\medskip
Here is an example that illustrates what the simple gluing theorem does. Suppose that $A$ and $B$ are the convex cores of two quasi-Fuchsian 3-manifolds which are immersed into some convex hyperbolic 3-manifold $M.$ By theorem (\ref{intersectgeomfinite}) the intersection of geometrically finite subgroups is geometrically finite  so $G=f_*\pi_1(A,a_0)\cap g_*\pi_1(B,b_0)$ is geometrically finite. By Scott's theorem this group is virtually embedded in $A$ and $B.$  

The induced gluing theorem applied with $G$ being this intersection produces finite covers of $A$ and $B$  where the pre-image of this intersection group is embedded. It then identifies the covers along submanifolds corresponding to these subgroups. We call these submanifolds (and the subsurfaces of the quasi-Fuchsian surfaces they correspond to) the {\em region of parallelism} between $A$ and $B.$  Assume for simplicity that $A$ and $B$ are embedded and contain spines which are disjoint surfaces $S_A,S_B.$ Then the region of parallelism corresponds to the maximal $I$-bundle between  $S_A$ and $S_B$ for which some fiber connects the basepoints of $S_A$ and $S_B.$ In general one must first pass to some finite cover since the region of parallelism does not correspond to a $\pi_1$-injective submanifold of $A$ and $B.$

To see this, since $A$ and $B$ are convex, every element of $G$ is represented by based geodesics in both $A$ and $B.$ The images of these geodesics under $f$ and $g$ are therefore equal. Hence the induced gluing will identify their lifts in $\tilde{A}$ and $\tilde{B}.$ 

It should be pointed out that the hyperbolic 3-manifold $S(\tilde{f},\tilde{g})$ produced by the simple gluing is not convex and in general the induced map $S(\tilde{f},\tilde{g})\rightarrow M$ is not $\pi_1$-injective. To achieve $\pi_1$-injectivity one needs some extra hypotheses, as in the convex combination theorem.

\section{\label{virtualtheorems}The virtual amalgam and virtual convex combination theorems.}
In this section we give two results which have the same conclusion as the convex combination theorem but with different hypotheses. The idea is that in some situations when one wishes to glue  two convex manifolds  the hypotheses of the convex combination theorem are always satisfied by certain finite covers of the manifolds  in question.

\medskip\noindent{\bf Definition.} Suppose that $M=M_1\cup M_2$ is a hyperbolic $n$-manifold which is the union of two convex hyperbolic $n$-submanifolds, $M_1,M_2.$ Suppose that $\pi:\tilde{M}\rightarrow M$ is a finite cover and $\tilde{M}_1$ is a component of $\pi^{-1}(M_i)$ and $\tilde{C}$ is a component of $\pi^{-1}(M_1\cap M_2).$  The hyperbolic $n$-manifold obtained from the disjoint union of $\tilde{M}_1$ and $\tilde{M}_2$ by identifying the copy of $\tilde{C}$ in each is called  a {\em virtual gluing} of $M_1$ and $M_2.$ This equals the simple gluing $S(\pi|\tilde{M}_1,\pi|\tilde{M}_2)$ with basepoints chosen in $\tilde{C}.$

\medskip
\begin{theorem}[virtual compact convex combination theorem]\label{virtualcombination}
Suppose that $M=M_1\cup M_2$ is a hyperbolic $n$-manifold, $M_1,M_2$ are compact, convex hyperbolic $n$-manifolds, and $C$ is a component of $M_1\cap M_2.$ Also suppose that $\pi_1C$ is a separable subgroup of both $\pi_1M_1$ and $\pi_1M_2.$ Then there is a virtual gluing, $N,$ of $M_1$ and $M_2$ along $C$ which has a convex thickening. In particular $N$ is isometric to a submanifold of ${\mathbb H}^n/hol(\pi_1N).$
\end{theorem}
\demo We will show that there is a virtual gluing of $M_1$ and $M_2$ which extends to a virtual gluing of the $\kappa$-thickenings. The result then follows from the convex combination theorem ({\ref{convexcombination}).

By (\ref{convexintersect}) $C$ is convex and thus $\pi_1C$ injects into $\pi_1M$ by (\ref{convexinjective}). Let $\pi:\tilde{M}\rightarrow M$ be the cover corresponding to $\pi_1C$ and let $\tilde{C}$ be a lift of $C$ to $\tilde{M}.$ Let $\tilde{M}_i$ be the component of $\pi^{-1}(M_i)$ which contains $\tilde{C},$ and set $N=\tilde{M}_1\cup\tilde{M}_2\subset\tilde{M}.$ Then $\pi_1N\cong \pi_1\tilde{C}\cong\pi_1C$ and so $N$ is isometric to a submanifold of ${\mathbb H}^n/\pi_1C.$
Define $$N^+ = \{\ x\in {\mathbb H}^n/\pi_1C\ :\ d(x,N)\le {\kappa}\ \}.$$
Then $N^+=T_{\kappa}(\tilde{M}_1) \cup T_{\kappa}(\tilde{M}_2).$ Clearly there is a covering map $T_{\kappa}(\tilde{M}_i)\rightarrow T_{\kappa}(M_i).$ 

Let  $\tilde{C}^+$ be the component of  $T_{\kappa}(\tilde{M}_1)\cap T_{\kappa}(\tilde{M}_2)$ which contains $C.$ Since $T_{\kappa}(\tilde{M}_1)$ and $T_{\kappa}(\tilde{M}_2)$ are convex, (\ref{convexintersect}) implies $\tilde{C}^+$ is convex. From (\ref{convexinjective}) it follows that $\pi_1\tilde{C}^+\rightarrow\pi_1({\mathbb H}^n/\pi_1C)$ is injective. Since $\tilde{C}\subset\tilde{C}^+$ it follows that  the inclusion $\tilde{C}\hookrightarrow\tilde{C}^+$ induces an isomorphism of fundamental groups. Thus $\tilde{C}^+$ is a convex thickening of $C.$    

We now show that $\tilde{C}^+$ is compact. Otherwise there is geodesic ray $\lambda$ in $\tilde{C}^+$ which starts at a point  $p\in\tilde{C}$ and leaves every compact set. Consider $\lambda_i\equiv \lambda\cap\tilde{M}_i.$ Since $p\in \tilde{C}\subset\tilde{M}_i$ it follows that $\lambda_i$ contains $p.$ By considering the universal cover of the $K$-thickening of $\tilde{M}_i,$ and using convexity of $\tilde{M}_i,$ it is easy to see that $\lambda_i=\lambda.$ Hence $\lambda\subset \tilde{C}.$ But this contradicts that $\tilde{C}$ is compact, and proves $\tilde{C}^+$ is compact.

{\bf Claim.}  There is a finite cover $Y_i\rightarrow T_{\kappa}(M_i)$ such that the natural map $p_i:\tilde{C}^+\rightarrow T_{\kappa}(M_i)$ lifts to an injective map
$\tilde{p}_i:\tilde{C}^+\rightarrow Y_i.$

Assuming the claim, consider the hyperbolic manifold $Y = Y_1\cup Y_2$ obtained by gluing $Y_1$ and $Y_2$ along $\tilde{C}^+.$ Then $Y$ contains $\tilde{M}=\tilde{M}_1\cup\tilde{M}_2$ and $Y_i=T_{\kappa}(\tilde{M}_1).$ The convex combination theorem now gives the result.

{\bf Proof of claim.} Choose a basepoint $x\in C$ and observe that $x\in C\subset M_i\subset T_{\kappa}(M_i).$ Under the identification given by the lift, $C\equiv\tilde{C},$ the point $x$ determines a basepoint $\tilde{x}\in\tilde{C}\subset\tilde{C}^+.$   In what follows all fundamental groups are based at the relevant basepoint.

  Let $S_i$ be the subset of $\pi_1T_{\kappa}(M_i)$ represented by loops of length at most $3$ times the diameter of $\tilde{C}^+.$  Since $\tilde{C}^+$ is compact  $S_i$ is finite. Since $\pi_1C$ is a separable subgroup of $\pi_1M_i$ there is finite index subgroup of  $\pi_1M_i$ which contains $\pi_1C$ and contains no element of $S_i\setminus\pi_1C.$ Let $Y_i\rightarrow T_{\kappa}(M_i)$ be the corresponding cover. 

Clearly the natural map $p_i$ lifts to $\tilde{p}_i:\tilde{C}^+\rightarrow Y_i.$ We have a basepoint $\tilde{p}_i(\tilde{x})\in Y_i.$ If $\tilde{p}_i$ is not injective there are distinct points $\tilde{x}_1\ne \tilde{x}_2\in\tilde{C}^+$ with $\tilde{p}_i(\tilde{x}_1)=\tilde{p}_i(\tilde{x}_2).$ Then there is a geodesic segment $\tilde{\alpha}$ in $\tilde{C}^+$ with endpoints $\tilde{x}_1$ and $\tilde{x}_2$, of length at most the diameter of $\tilde{C}^+.$ Observe that $\tilde{\alpha}$ maps to a loop $\tilde{p}_i\circ\tilde{\alpha}$ in $Y_i.$ Thus $\tilde{p}_i\circ\tilde{\alpha}$ projects to a loop, $\alpha,$ in $T_{\kappa}(M_i)$ based at the point $x_1 = p_i(\tilde{x}_1).$ Let $\tilde{\beta}$ be a geodesic segment of minimal length in $\tilde{C}^+$ with endpoints $\tilde{x}$ and $\tilde{x}_1.$ Then $\beta = p_i\circ\tilde{\beta}$ is a geodesic segment in $T_{\kappa}(M_i)$  starting at $x$ and ending at $x_1.$  Thus $\gamma = \beta.\alpha.\beta^{-1}$ is a loop in $T_{\kappa}(M_i)$ based at $x$ of length at most $3$ times the diameter of $\tilde{C}^+.$  Thus $[\gamma]$ is an element of $S_i.$  It is easy to see that $\gamma$ lifts to a loop in $Y_i$ based at $\tilde{p}_i(\tilde{x})$ so  $[\gamma]\in Im[\pi_1C\rightarrow\pi_1T_{\kappa}(M_i)].$ This implies $\gamma$ lifts to a loop in $\tilde{C}^+$ based at $\tilde{x}$ which contradicts that $\tilde{\alpha}$ has distinct endpoints. This proves the claim and the theorem.\qed

\begin{proposition}[increasing the rank of a cusp]\label{increaserank} Suppose that $M$ is a convex hyperbolic $n$-manifold and $C$ is a thin cusp contained in $M.$ Suppose that   $\Gamma< Isom({\mathbb H}^n)$ is a discrete group of parabolic isometries which contains $hol(\pi_1C).$  Then there is a finite-index subgroup $\Gamma'<\Gamma$ which contains $hol(\pi_1C)$ and a horoball $D\subset{\mathbb H}^n$ stabilized by $\Gamma'$ with the following property. Set $Q=D/\Gamma'$ then there is a hyperbolic $n$-manifold $N=M\cup Q$ such that $M\cap Q=C.$ Furthermore $N$ has a convex thickening.\end{proposition}
\demo  Let $\pi:\tilde{M}\rightarrow M$ be the universal cover. Use the developing map to isometrically identify $\tilde{M}$ with a subset of ${\mathbb H}^n.$ Let $\tilde{C}$ be a component of $\pi^{-1}C.$ Let $D$ be the horoball which contains $\tilde{C}$ and so that $\partial D$ contains $\tilde{C}\cap\pi^{-1}(\partial_cC).$
Let $D^-\subset D$ be the smaller horoball such that $D$ is a $\kappa$-neighborhood of $D^-.$ Then
$C^-=(\tilde{C} \cap D^-)/hol(\pi_1C)$ is a smaller cusp contained in $C.$ The cusp $W=D/\Gamma$ has boundary $\partial W$ which is a Euclidean manifold. There is also a smaller cusp $W^-=D^-/\Gamma.$

Given $K>0$ let $Y_1=T^{rel}_{K}(M;C^-).$ We regard the universal cover $\tilde{Y}_1$ as a subset of ${\mathbb H}^n$ so that it contains $\tilde{M}$ in the natural way. Let $C_1$ be the image of $D\cap\tilde{Y}_1$  under projection to $Y_1.$ This is a cusp in $Y_1$ and by choosing $K$ sufficiently large (and $K\ge\kappa$) we may arrange that $C_1$ contains $C$ and so $\partial_cC\subset\partial_cC_1.$ Thus $C_1$ is the cusp  in $Y_1$ which naturally corresponds to the cusp $C$ in $M_1.$  Since $hol(\pi_1C)<\Gamma$ there is a natural local-isometry $f:C_1\rightarrow W.$

\medskip{\bf Claim.} There is a finite cover $\tilde{W}\rightarrow W$ so that $f$ lifts to an embedding $\tilde{f}:C_1\rightarrow\tilde{W}.$

\medskip

Assuming this claim, since $C^-\subset C\subset C_1$ it follows that $\tilde{f}| : C^-\rightarrow \tilde{W}^-$ is  injective. In order to fit with the notation used in the convex combination theorem we now use $M_1$ to denote $M$ and $M_2$ to denote $\tilde{W}^-.$  We use $\tilde{f}$ to identify $C^-\subset M_1$ with its image in $M_2$ and set $M=M_1\cup M_2$ thus $M_1\cap M_2 = C^-.$    Similarly let $Y_2=T_{\kappa}(M_2)=\tilde{W}$ and use $\tilde{f}$ to identify $C_1\subset Y_1$ with its image in $Y_2$ then $Y=Y_1\cup Y_2$ and $Y_1\cap Y_2 = C_1.$

\centerline{\epsfysize=40mm
\epsfbox{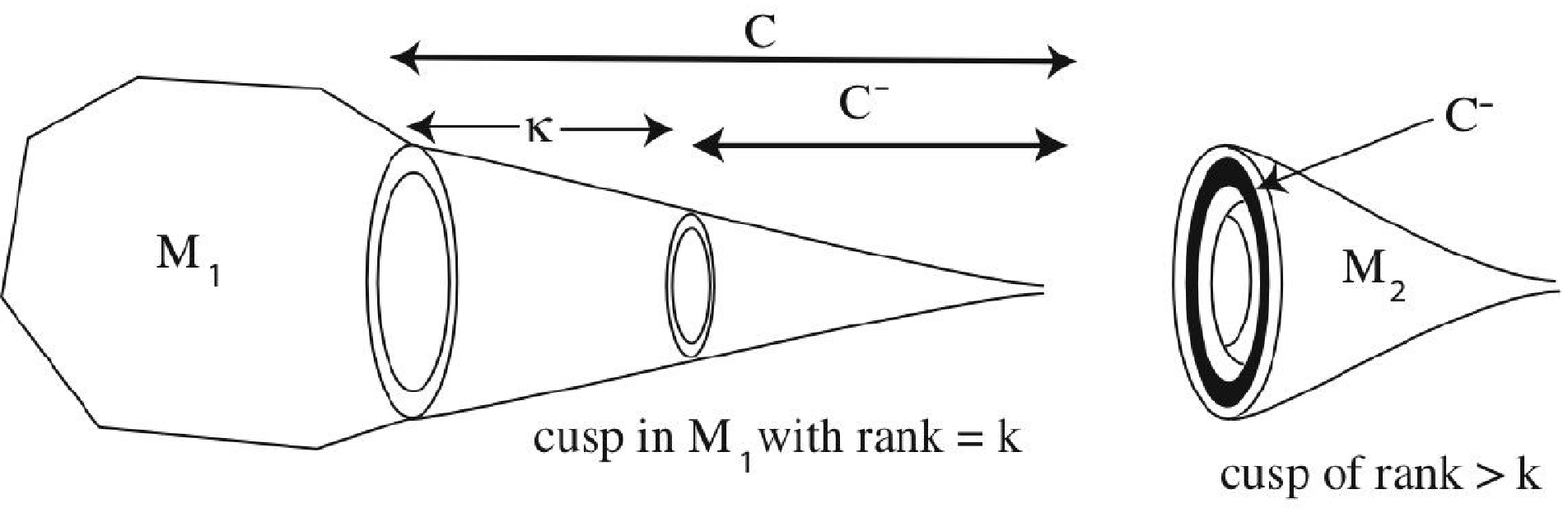}}   
   
We now check the 6 hypotheses of the convex combination theorem are satisfied.  Using the fact that a relative thickening of a convex manifold is convex, it is easy to check  that $M_1,M_2,Y_1,Y_2$ are all convex which gives conditions (1) and (2). Clearly $M$ is a submanifold of $Y$ and $Y_i$ is a thickening of $M_i$ which implies condition (3).   Since $K\ge\kappa$ by  (\ref{Trel})(c) if $x\in \overline{M_1\setminus M_2}$ then $exp_x:T_xM_1\rightarrow Y_1$ is defined on vectors of length at most $\kappa.$ Since $Y_2=T_{\kappa}(M_2)$ it follows that $exp_x:T_xM_2\rightarrow Y_2$ is defined on vectors of length at most $\kappa.$ This implies conditions (4) and  (5) are  satisfied. Finally, $C^-=M_1\cap M_2$ is contained in $C_1=Y_1\cap Y_2,$ and both are connected, so condition (6) is satisfied. The convex combination theorem implies that $M$ has a convex thickening $T_{\infty}(M).$ This contains $Y$ and thus contains $M\cup Y_2.$ Now $M\cap Y_2=C$ and this proves the theorem with $Q=Y_2.$

\medskip
{\bf Proof of Claim.}  Since $C$ is a thin cusp $C_1$ is also a thin cusp. By (\ref{Trel})(e) there is a product cusp, $P,$ which contains $C_1.$ Clearly the local isometry $f$ has an extension to a local isometry $f:P\rightarrow W.$ Thus it suffices to prove this
extension is injective. A product cusp has a $1$-dimensional foliation by rays orthogonal to the cusp boundary. If two distinct points in $P$ have the same image under $f$ then the rays through these points have the same image under $f.$ It follows that if $f$ is not injective there are two points in the cusp boundary $\partial_cP$ which have the same image under $f.$ Thus it suffices to show there is a finite cover of $\partial W$ so that $f|:\partial_cP\rightarrow\partial W$ lifts to an embedding. 

Since $W$ is a cusp  $\pi_1W\cong\pi_1\partial W.$  The manifold $\partial W$ is Euclidean and so has a fundamental group which is virtually free-abelian. Thus $\pi_1W$ is subgroup separable. Since $P$ is a product cusp $\partial_cP$ is compact.  Since $f$ is a local isometry and $\partial_cP$ is compact and convex it follows from a standard argument that such a cover exists.
\qed

\medskip

The convex combination theorem sometimes enables one to glue geometrically finite manifolds together to obtain a geometrically finite manifold. This corresponds to forming an amalgamated free product of two geometrically finite groups, amalgamated along their intersection.

\medskip\noindent{\bf Definition.} Two subgroups $A,B$ of a group $G$ can be {\em virtually amalgamated} if there are finite index subgroups $A'<A$ and $B'<B$ such that the subgroup, $G',$ of $G$ generated by $A'$ and $B'$ is the free product of $A'$ and $B'$ amalgamated along $A'\cap B'.$ We also say that $G'$ is a {\em virtual amalgam} of $A$ and $B.$

\medskip\noindent{\bf Definition.} Two non-trivial parabolic subgroups $\Gamma_1,\Gamma_2 < Isom({\mathbb H}^n)$ are called {\em compatible} if either:\\
(1) $\Gamma_1$ and $\Gamma_2$ stabilize distinct points on the sphere at infinity, or\\
(2) $\Gamma_1\cap\Gamma_2$ has finite index in at least one of the groups $\Gamma_1$ or $\Gamma_2.$ 

\medskip
The second condition is equivalent to saying that, up to taking subgroups of finite index, that one group is a subgroup of the other. Two discrete groups $\Gamma,\Gamma'<Isom({\mathbb H}^n)$ have {\em compatible parabolic subgroups} if every maximal parabolic subgroup of $\Gamma$ is compatible with every maximal parabolic subgroup of $\Gamma'.$

\begin{theorem}[GF subgroups have virtual amalgams]\label{virtualamalgam}  Suppose that $\Gamma$ is a discrete subgroup of $Isom({\mathbb H}^n).$ Suppose that $\Gamma_1$ and $\Gamma_2$ are two geometrically finite subgroups of $\Gamma$ which have compatible parabolic subgroups. Then $\Gamma_1$ and $\Gamma_2$ can be virtually amalgamated and the result is a geometrically finite group.
\end{theorem}
\demo Since $\Gamma$ is linear there is a torsion-free subgroup of finite index in $\Gamma.$ We may replace each of $\Gamma_1,\Gamma_2,\Gamma$ by their intersections with this subgroup. Thus we may assume that $\Gamma$ is torsion-free.

Let $N={\mathbb H}^n/\Gamma.$  Choose a basepoint $\tilde{x}\in {\mathbb H}^n$ and let $x$ be the image in $N$ of $\tilde{x}.$ This choice determines an identification  $\pi_1(N,x) \equiv \Gamma.$ Let $H_i = CH(\Gamma_i\cdot \tilde{x})$ then
$M_i = H_i/\Gamma_i$  is a convex hyperbolic manifold. Let $x_i$ be the image of $\tilde{x}$ in $M_i.$ The choice of $\tilde{x}$ determines an identification $\pi_1(M_i,x_i) \equiv \Gamma_i.$ The inclusion $H_i\subset{\mathbb H}^n$ covers a  local isometry $\rho_i:M_i\rightarrow N.$ Then $\rho_{i*}$ maps $\pi_1(M_i,x_i)$ into a subgroup of $\pi_1(N,x)$ and under the identifications $\pi_1(N,x) \equiv \Gamma$ and $\pi_1(M_i,x_i) \equiv \Gamma_i$ the map $\rho_{i*}$ is inclusion $\Gamma_i\subset \Gamma.$

\medskip
{\bf Case 1.} $M_1$ and $M_2$ have no cusps.

\medskip  Define $Y_i=T_{\kappa}(M_i).$  We apply the simple gluing theorem to the maps $\rho_i:(Y_i,x_i)\rightarrow (N,x).$ Thus there are finite covers $\tilde{Y}_i\rightarrow Y_i$ and a simple gluing $\tilde{Y}=\tilde{Y}_1\cup\tilde{Y}_2.$ In addition there are lifts of the basepoints, $\tilde{x}_i,$ to $\tilde{Y}_i$ and these lifted basepoints are identified.  By (\ref{simplegluelemma})(e) $\tilde{Y}_1\cap\tilde{Y}_2$ is connected. Obviously it contains the basepoint. Let $\tilde{M}_i\subset\tilde{Y}_i$ be the corresponding cover of $M_i.$ These also contain the basepoint thus $\tilde{M}_1\cap\tilde{M}_2\ne\phi.$ 

Since $\tilde{Y}_i=T_{\kappa}(\tilde{M}_i)$ conditions (1)-(5) of the convex combination theorem are satisfied. The above remarks imply condition (6) of the convex combination theorem is satisfied.  Since $\Gamma_i$ is geometrically finite, $M_i$ has finite volume thus so does $\tilde{M}_i.$ The convex combination theorem then implies that $\tilde{M}_1\cup \tilde{M}_2$ has a convex thickening which is geometrically finite.  Also (\ref{freeproduct}) implies that $\pi_1(\tilde{M}_1\cup \tilde{M}_2)$ is an amalgamated free product. This completes the proof of case 1.

\medskip
The same proof works unchanged when every cusp of $N_1$ and of $N_2$ has maximal rank.  This is because thickening and relative thickening give the same result if and only if all the cusps have maximal rank.  However when some cusp has rank less than maximal, we must change the proof to use relative thickening instead of thickening in order to apply the virtual simple gluing theorem. This leads to some technicalities which we now address.

\medskip
{\bf Case 2.} For every cusp, $C_i,$ of $M_i,$ the cusp, $C,$ in $N$ which contains the image of $C_i$ has the same rank as $C_i.$

\medskip
Let ${\cal C}$ be a maximal collection of disjoint cusps in $N.$ Then ${\cal C}_i = \rho_i^{-1}({\cal C})$ is a maximal collection of disjoint cusps in $M_i.$  Let $C_{i,j}$ be a component of ${\cal C}_{i}.$ Then $\partial_cC_{i,j}$ is compact because $M_i$ has finite volume. Given $K\ge\kappa$ 
define $Y_i=T_{K}^{rel}(M_i;{\cal C}_i).$ There is  a local isometry $\rho_i:Y_i\rightarrow N$ which extends $\rho_i:M_i\rightarrow N.$ We now choose $K$ so large that the following cusp condition is satisfied for all $i \in \{1,2\}$ and all $j.$ Suppose that $C^*_{i,j}$ is a cusp in $Y_i$ corresponding to a cusp $C_{i,j}\in{\cal C}_i.$ If the images of $C^*_{1,j}$ and $C^*_{2,j}$ are contained in the same cusp $C\in{\cal C}$ of $N$ then $\rho_1(\partial_c C_{1,j})\subset \rho_2(C^*_{2,j})$ and $\rho_2(\partial_c C_{2,j})\subset \rho_1(C^*_{1,j}).$

We then obtain finite covers $p_i:\tilde{Y}_i\rightarrow Y_i$ and $\tilde{M}_i=p_i^{-1}(M_i)$  and simple gluings $\tilde{Y}=\tilde{Y}_1\cup\tilde{Y}_2$ and $\tilde{M}=\tilde{M}_1\cup\tilde{M}_2$ as in case (1).  Also $\tilde{Y}_1\cap\tilde{Y}_2$ and $\tilde{M}_2\cap\tilde{M}_2$ are connected by (\ref{simplegluelemma})(e) and convex by (\ref{convexintersect}) and contain the basepoint. However since we only used relative thickenings condition (4) of the convex combination theorem is not satisfied.  Define $Y_i^+=T_{K}(M)$ thus $Y_i\subset Y_i^+.$ Let $\tilde{Y}_i^+$ be the corresponding covering of $Y_i^+.$ 

\medskip
{\bf Claim 1.}   The simple gluing of $\tilde{Y}= \tilde{Y}_1\cup\tilde{Y}_2$ extends to a simple gluing  $\tilde{Y}^+= \tilde{Y}_1^+\cup\tilde{Y}_2^+.$

\medskip
\noindent Assuming the claim the hypotheses of the convex combination theorem are satisfied by $\tilde{M}=\tilde{M}_1\cup\tilde{M}_2$ and $\tilde{Y}^+=\tilde{Y}_1^+\cup\tilde{Y}_2^+.$ The remainder of the proof of case 2 is the same as case 1.\medskip

\noindent{\bf Proof of claim 1.}  There is a local isometry $f:\tilde{Y}\rightarrow N.$  Let $C\in{\cal C}$ be a cusp of $N$ and $A\subset\tilde{Y}$ a component of $f^{-1}(C).$ Suppose that $D_i$ is a component of $A\cap\tilde{Y}_i.$ Then $D_i$ is a cusp in $\tilde{Y}_i.$

\medskip
{\bf Claim 2.}  There is a finite cover $\tilde{C}\rightarrow C$ such that the map $f|:A\rightarrow C$ lifts to an embedding $\tilde{f}:A\rightarrow\tilde{C}.$

\medskip
\noindent{\bf Proof of claim 2.} If $A\cap\tilde{Y}_2$ is empty then $A = D_1\subset\tilde{Y}_1$ is a cusp in $\tilde{Y}$ and in particular is convex. By hypothesis  $f_*\pi_1(D_1)$ has finite index in $\pi_1(C)$ so the cover $\tilde{C}\rightarrow C$ corresponding to this subgroup is finite and $f$ lifts. Since $D_1$ is convex $\tilde{f}:D_1\rightarrow\tilde{C}$ is injective. A similar argument works if $A\cap\tilde{Y}_1$ is empty.

So we may suppose $D_1$ and $D_2$ are components as above and $D_1\cap D_2\ne\phi.$ Choose $x\in D_1\cap D_2$ then $G = f_*\pi_1(D_1,x)\cap f_*\pi_1(D_2,x)$ is the intersection of finite index subgroups thus has finite index in $\pi_1(C,fx).$
Let $\tilde{C}\rightarrow C$ be the cover corresponding to $G.$ Then there are induced finite covers $\tilde{D}_i\rightarrow D_i$ and maps $\tilde{f}_i:\tilde{D}_i\rightarrow \tilde{C}$ covering $f|D_i$ which are  injective (for the same reason as before).

Let $Z = S(\tilde{f}_1,\tilde{f}_2)$ be the induced gluing of $\tilde{D}_1$ and $\tilde{D}_2$ using the copy of the basepoint $x$ in both $D_1$ and $D_2.$ Since $\tilde{f}_1,\tilde{f}_2$ are both injective this is a simple gluing, so $Z = \tilde{D}_1 \cup \tilde{D}_2$ identified along the pre-images of the component of $\tilde{f}_1(\tilde{D}_1)\cap\tilde{f}_2(\tilde{D}_2)$ which contains the image of the basepoint.

\medskip{\bf Claim 3.} $\tilde{f}_1(\tilde{D}_1)\cap\tilde{f}_2(\tilde{D}_2)$ is connected. 

\medskip\noindent{\bf Proof of claim 3.} Otherwise there is an essential loop, $\gamma$ in $\tilde{f}_1(\tilde{D}_1)\cup\tilde{f}_2(\tilde{D}_2)$ which is trivial in $\tilde{C}.$ To see this, choose $\gamma=\gamma_1\cdot\gamma_2$ so that $\gamma_i\subset\tilde{f}_i(\tilde{D}_i)$ and the endpoints of $\gamma_i$ are in different components of $\tilde{f}_1(\tilde{D}_1)\cap\tilde{f}_2(\tilde{D}_2).$ Now using that $\tilde{f}_{1*}\pi_1(\tilde{D}_1)\cong\pi_1(C)$ we can compose $\gamma$ with a loop in $\tilde{f}_{1*}\pi_1(\tilde{D}_1)$ so that it is contractible in $\tilde{C}.$

In the universal cover, $X,$ of $\tilde{C}$ the loop $\gamma$ is covered by a loop $\tilde{\gamma} = \tilde{\gamma}_1\cdot\tilde{\gamma}_2.$ There are copies, $E_i,$ of the copy of the universal cover of $\tilde{f}_i(\tilde{D}_i)$ in $X,$ such that $\tilde{\gamma}_i\subset E_i.$ Thus $E_1\cap E_2$ is not connected. But they are convex sets in $X\subset{\mathbb H}^n$ which is impossible. This proves claim 3.\qed

\medskip
\noindent{\bf Proof of claim 2 resumed.}  By claim 3 we can identify $Z$ with $\tilde{f}_1(\tilde{D}_1)\cup\tilde{f}_2(\tilde{D}_2) \subset\tilde{C}.$   There is a natural map $\theta:Z\rightarrow A$ which is surjective.
The cusp $p_1(D_1)$ of $Y_1$ contains a unique  cusp $C_1$ of $M_1$ where $C_1\in{\cal C}_.$ Let $E$ be the pre-image of $\partial_c C_1$ in $\tilde{D}_1.$
 Now $C_1$ deformation retracts to $\partial_c C_1,$  and $\tilde{D}_1$ is a covering of a thickening of $C_1,$ thus $\tilde{D}_1$ deformation retracts to $E.$ Hence $E$ is connected. By choice of $K$ we have $\rho_1(\partial_c C_1)\subset \rho_2(p_2(D_2)).$ It follows that $E$ equals the pre-image in $\tilde{D}_2$ of $\rho_2^{-1}(\rho_1\partial_c C_2).$ Thus   $E\subset \tilde{D}_1\cap \tilde{D}_2,$ and so $ \theta(E) \subset D_1 \cap D_2.$
 
Since $D_1$ deformation retracts to $E_1,$ it follows that  $incl_*:\pi_1E_1\cong\pi_1D_1.$ Since $D_i$ is convex, $f|:D_i\rightarrow C$ is $\pi_1$-injective, so we may regard $\pi_1D_i$ as a subgroup of $\pi_1C.$ Now $E_1$ carries $\pi_1D_1$ and $E_1\subset D_2$ thus $\pi_1D_1$ is a subgroup of $\pi_1D_2.$ The reverse inclusion is obtained similarly. Hence $\pi_1D_1\cong\pi_1D_2.$ 
The above argument now applies with $G = \pi_1D_1\cap\pi_1D_2 = \pi_1D_i.$ Then $\tilde{C}=C$ and $Z = f_1(D_1) \cup f_2(D_2) = A \subset C.$ This proves claim 2.\qed

\medskip  We can now glue $\tilde{C}$ onto $\tilde{Y}$ using $g$ to obtain $\tilde{Y}\cup\tilde{C}$ with $\tilde{Y}\cap\tilde{C}=A.$ We do this to for every cusp $C$ of $N$ and every component $A$ of $f^{-1}(C).$ The result is a manifold $Z$ obtained from $Y$ by adding a complete cusp onto every cusp of $Y.$  Clearly $\tilde{Y}^+_i=T_{K}(\tilde{M}_i)$ is contained in $Z,$ thus $\tilde{Y}^+=\tilde{Y}_1^+\cup\tilde{Y}_2^+\subset Z.$ This proves claim 1 and complete the proof of case 2.\qed

\medskip
{\bf General Case.}

\medskip

Given a cusp $C_1\in{\cal C}_1$ of $M_1$ let $C\in{\cal C}$ be the corresponding cusp of $N$ which contains the image of $C_1.$ If $rank(C_1) < rank(C)$ we can use (\ref{increaserank}) to glue a finite cover of $C$ onto $C_1$ and thicken to obtain a convex manifold.  After doing this to every cusp of $M_1$ we obtain  a finite volume convex manifold $M_1^+.$ Similarly we obtain $M_2^+.$ Using case 2 there are finite covers $\tilde{M}_1^+,\tilde{M}_2^+$ and a simple gluing $\tilde{M}^+=\tilde{M}_1^+\cup\tilde{M}_2^+$ which has a convex thickening, $P,$ of finite volume.

The final step is to take a certain submanifold of a certain covering space of $P$ to get the required manifold. It is only  at this last step that the hypothesis of compatible parabolic subgroups is used.

The manifold $\tilde{M}_i^+$ contains a submanifold $\tilde{M}_i$ which is a finite cover of $M_i.$ Set $W=\tilde{M}_1\cup\tilde{M}_2\subset\tilde{M}^+.$ Suppose $\tilde{C}_i$  is a cusp in $\tilde{M}_i$ for some $C_i\in{\cal C}_i.$ The compatibility hypothesis implies that if $\tilde{C}_1\cap \tilde{C}_2$ is not empty then $\tilde{C}_1\subset \tilde{C}_2$ or $\tilde{C}_2\subset \tilde{C}_1.$ Hence $\tilde{M}^+ = W \cup \bigcup_i D_i$ where $D_i$ is a cusp (which is isometric to a finite cover of some cusp, $C\in{\cal C}$ of $N$) glued onto $W$ along a cusp $\tilde{C}_1$ of $\tilde{M}_1$ or a cusp $\tilde{C}_2$ of $\tilde{M}_2.$
It follows that $W$ is $\pi_1$-injective in $\tilde{M}^+.$ 

We have $W\subset\tilde{M}^+\subset P.$  Let $\tilde{P}\rightarrow P$ be the cover corresponding to $\pi_1W.$ Let $\tilde{W}$ be the lift of $W$ to $\tilde{P}.$ Let $f:\tilde{P}\rightarrow N$ be the natural map. Let $X\subset\tilde{P}$ be the component of $f^{-1}(N\setminus int({\cal C}))$ which contains $\tilde{W}\setminus f^{-1}(int({\cal C})).$ Every point of $X$ is within a distance $\kappa$ of $\tilde{W}\setminus f^{-1}(int({\cal C})),$ therefore $X$ is compact. 

Each component of $X\cap f^{-1}({\cal C})$ is a compact set contained in the cusp boundary of some cusp of $\tilde{P}.$ Thus $CH(X) = X \cup \bigcup F_i$ where $F_i$ is a thin cusp. Now $CH(X)$ is convex, finite volume, and has a finite volume $\epsilon$-thickening, so it is a geometrically finite manifold. Also $W$ is contained in $CH(X).$ Thus $CH(X)$ is the desired manifold, which completes the proof.\qed

\section{Some Constructions of Hyperbolic Manifolds.}
In this section we use the convex combination theorem to give  constructions of geometrically finite hyperbolic manifolds in dimensions bigger than $3.$ The basic idea is to take two hyperbolic manifolds of dimensions $m,n$ each of which contains a copy of the same totally geodesic submanifold of dimension $p$ and then glue the manifolds along the submanifold and thicken to get a convex hyperbolic manifold of dimension $m+n-p.$

Consider a  hyperbolic $m$-manifold, $M={\mathbb H}^m/G,$ where $G$ is a discrete subgroup of $Isom({\mathbb H}^m).$ We first describe a way to thicken $M$ to obtain a hyperbolic $n$-manifold with $n>m.$ We may embed ${\mathbb H}^m$ isometrically as a totally geodesic subspace of ${\mathbb H}^n.$ We regard $Isom({\mathbb H}^m)$ (and thus $G$)  as a subgroup of $Isom({\mathbb H}^n).$ Given $R>0$ let $$N = \{\ x\in{\mathbb H}^n\ :\ d(x,{\mathbb H}^m)\le R\ \}.$$ Define $T(M;n,R)=N/G.$ This is a convex hyperbolic $n$ manifold with strictly convex boundary. There is a projection $\pi:{\mathbb H}^n\rightarrow{\mathbb H}^m$ given by the nearest point retraction. This map is $G$-equivariant thus we get a map $p:N/G\rightarrow M$ which is a Riemannian submersion and is a disc bundle over $M.$

\medskip
\noindent{\bf Example 1.}\\
 For $i\in\{1,2\}$ suppose  that $M_i$ is a closed hyperbolic 3-manifold which contains a simple closed geodesic $C_i$ of  length $\ell$ such that  the holonomy around $C_i$ is a pure translation. We suppose that $C_i$ has  a tubular neighborhood $V_i$ in  $M_i$ of radius $\kappa,$ the thickening constant. 

Consider the hyperbolic $5$-manifolds $M_i^+=T(M_i;5,\kappa)=N_i/G_i$ where $N_i$ is a $\kappa$-neighborhood of a hyperbolic $3$-space ${\mathbb H}^3_i\equiv\tilde{M}_i$ in ${\mathbb H}^5.$  We may choose ${\mathbb H}^3_1$  to be orthogonal to ${\mathbb H}^3_2$ and intersect along a geodesic which covers $C_1$ in $M_1$ and $C_2$ in $M_2.$ Then $\tilde{V}=N_1\cap N_2$ is a neighborhood of ${\mathbb H}^3_1\cap {\mathbb H}^3_2.$ Let $V_i\cong S^1\times D^4$ denote the image of $\tilde{V}$ in $M_i^+.$ Let $M^+$ be the non-convex hyperbolic $5$-manifold obtained by gluing $M_1^+$ to $M_2^+$ by identifying $V_1$ with $V_2$ in a way which is covered by the identifications between $N_1$ and $N_2.$
By the convex combination theorem $M^+$ has a convex thickening.

\medskip
\noindent{\bf Example 2.}\\
For $i\in\{1,2\}$  suppose that $M_i$ is a closed hyperbolic 3-manifold which contains a totally geodesic surface $F_i$ and assume that $F_1$ is isometric to $F_2.$ Also assume that $M_i$ contains a tubular neighborhood of $F_i$ of radius the thickening constant $\kappa.$ Now consider the hyperbolic 4-manifolds $M_i^+=T(M_i;4,\kappa).$ These may be glued by isometrically identifying neighborhoods of $F_1$ and $F_2$ to obtain a hyperbolic 4-manifold, $M^+=M^+_1\cup M^+_2,$ which is homotopy equivalent to the result of gluing $M_1$ to $M_2$ by identifying $F_1$ with $F_2.$ By the convex combination theorem $M^+$ has a convex thickening.

\medskip
\noindent{\bf Example 3.}\\
This time we will construct a hyperbolic 4-manifold by gluing a hyperbolic surface with geodesic boundary component, $C_1,$  along $C_1$ to a geodesic, $C_2,$ of the same length in a closed hyperbolic $3$-manifold and thickening. In order that the resulting $4$-manifold have a convex thickening it suffices that $C_1$ and $C_2$ both have tubular neighborhoods of radius $\kappa$ in their respective manifolds. In addition $C_2$ should be a pure translation.

This example may be modified to allow the holonomy along $C_2$ to have a small non-zero rotational part, by deforming the holonomy of the surface from a subgroup of $Isom({\mathbb H}^2)$ into a nearby subgroup of $Isom({\mathbb H}^4)$ so that $C_1$ and $C_2$ have the same holonomy.

\begin{theorem}[gluing hyperbolic manifolds along isometric submanifolds]\ \\
For $i\in\{1,2\}$ suppose that $M_i$ is a convex hyperbolic $m_i$-manifold. Suppose that $p<\min(m_1,m_2)$ and $P_i$ is a closed hyperbolic $p$-manifold in $M_i.$ Suppose that $f:P_1\rightarrow P_2$ is an isometry. Suppose that $P_i$ has a tubular neighborhood in $M_i$ of radius at least the thickening constant $\kappa.$ Suppose that the holonomy of $P_i$ is contained in the subgroup $Isom({\mathbb H}^p)$ of $Isom({\mathbb H}^{m_i}).$ Let $M$ be the space obtained by gluing $M_1$ to $M_2$ by using $f$ to identify $P_1$ with $P_2.$ Then $M$ has a convex thickening which is a hyperbolic $(m_1+m_2-p)$-manifold.
\end{theorem}

\begin{theorem}[virtual gluing hyperbolic manifolds along isometric submanifolds]\ \\
For $i\in\{1,2\}$ suppose that $M_i$ is a convex hyperbolic $m_i$-manifold. Suppose that $p<\min(m_1,m_2)$ and $P_i$ is a closed hyperbolic $p$-manifold in $M_i.$ Suppose that $f:P_1\rightarrow P_2$ is an isometry. Suppose that the holonomy of $P_i$ is contained in the subgroup $Isom({\mathbb H}^p)$ of $Isom({\mathbb H}^{m_i}).$ Then there are finite covers $\tilde{M}_i$ of $M_i$ and lifts $\tilde{P}_i$ of $P_i$ to $\tilde{M}_i$ with the following property. Let $\tilde{M}$ be the space obtained by gluing $\tilde{M}_1$ to $\tilde{M}_2$ by using $f$ to identify $\tilde{P}_1$ with $\tilde{P}_2.$ Then $\tilde{M}$ has a convex thickening which is a hyperbolic $(m_1+m_2-p)$-manifold.
\end{theorem}

\demo By (\ref{totgeodesic}) $\pi_1P_i$ is a separable subgroup of $\pi_1M_i.$ It follows  there is a finite cover, $\tilde{M}_i$ of $M_i$ and a lift of $P_i$ to $\tilde{M}_i$ which has a tubular neighborhood of radius $\kappa.$ The result follows from the previous theorem.\qed 

\section{Subgroup Separability.} 

Two groups $G_1$ and $G_2$ are {\em commensurable} if there is a group $H$ which is isomorphic to finite index subgroups of both $G_1$ and $G_2.$ Two path-connected topological spaces $X,Y$ are {\em commensurable} if there are finite sheeted covers $\tilde{X},\tilde{Y}$ which are homeomorphic. Clearly commensurable spaces have commensurable fundamental groups.

A subgroup $H$ of a group $G$ is {\em separable in $G$}  if for every $g\in G\setminus H$
there is a subgroup of finite index $K<G$ such that $H\le K$ and $g\notin K.$  The group $G$ is {\em subgroup separable} if every subgroup is separable and is {\em LERF} if every finitely generated subgroup is separable. 

\medskip
{\bf Definition.} Suppose that $F$ is a compact, connected, surface with $\chi(F)<0.$ Let $\{\partial_i F\}_{1\le i\le n}$ denote the boundary components of $F$ and let $\{T_i\}_{1\le i\le n}$ denote a collection of distinct tori.  The {\em tubed surface}, $X,$ obtained from $F$  is the $2$-complex $X$ obtained obtained by homeomorphically identifying each component $\partial_iF$  with an essential simple closed curve on the torus $T_i.$ 

\medskip
Observe that
$\pi_1X$ is a topological realization of a graph of groups with cyclic edge
groups and with vertex groups that are either
${\mathbb Z}^2$ or finitely generated free groups. It follows that $X$ is a
$K(\pi,1).$
\begin{lemma}[tubed surface is LERF]\label{tubedisLERF} If $X$ is a tubed surface then $\pi_1X$ is
LERF.\end{lemma}
\demo We first show that all tubed surfaces are commensurable. Let $A$ denote the compact surface obtained by deleting the interior of a disc from a torus. It is easy to check that up to commensurability  every  compact, connected, surface with negative Euler characteristic and non-empty boundary is commensurable with $A.$ Let $Y$ denote the tubed surface obtained from $A$.  It follows that every tubed surface is commensurable with with $Y.$ Thus all tubed surfaces have commensurable fundamental groups. 

R. Gitik proved in theorem (4.4) of  \cite{G1} that an amalgam of a free group, $F,$ and a LERF group, $H,$ is again LERF, provided the amalgamating subgroup is a maximal cyclic subgroup of  $F$.  It is easy to see that ${\mathbb Z}\oplus{\mathbb Z}$ is LERF and that $\pi_1\partial A$ is a maximal cyclic subgroup of $\pi_1A.$ It follows that $\pi_1Y$ is LERF.  Scott shows in  lemma (1.1) of \cite{SC1} that the property of being LERF is a commensurability invariant and it follows that the fundamental group of every tubed surface is LERF.
\qed

\medskip
{\bf Definiton.} If $X$ is a path connected topological space and $G$ is a subgroup of $\pi_1(X)$
 we say that {\it $G$ is virtually embedded in $X$} if there is a finite
sheeted cover $p:\tilde{X}\rightarrow X$ and a
path-connected subspace
$Y\subset\tilde{X}$ such that $p|_*:\pi_1Y \rightarrow \pi_1X$ is injective with image $G.$

\medskip
In \cite{SC1} Scott introduced the notion of virtually embedded for surface subgroups (but he used the term {\em almost geometric} which, in the context of 3-manifolds, has certain connotations we prefer to avoid) and
used subgroup separability to prove that finitely generated subgroups of surface
groups are virtually embedded.  Scott's result extends to separable subgroups of  3-manifolds:
\begin{theorem}[virtually embedded subgroups]\label{virtgeom} Suppose that $M$ is a connected 3-manifold and $G<\pi_1M$ is a
finitely generated separable subgroup. Then $G$ is a  virtually embedded
subgroup.
\end{theorem}
\demo The following argument is standard. Let $p_G:\tilde{M}_G\rightarrow M$ be the cover corresponding to $G.$
By the compact core theorem, \cite{HE2}, there is a compact submanifold $Y\subset
\tilde{M}_G$ such that the inclusion, $\iota:Y\rightarrow\tilde{M}_G,$ is a
homotopy equivalence. We may choose a triangulation of
$M$ such that $Y$ is a finite subcomplex of the induced triangulation of
$\tilde{M}_G.$ We will show how to construct a tower of finite covers (as in the proof of the
loop theorem)
$$p_n:\tilde{M}_{n+1}\rightarrow\tilde{M}_n$$ 
with $\tilde{M}_0=M$ and such that the map $f_0=p_G\circ\iota:Y\rightarrow M$ lifts
to each cover $f_n:Y\rightarrow\tilde{M}_n.$ The singular set of $f_n$ is
$$S(f_n)\ =\ \{\ x\in Y\ :\ \#|f_n^{-1}(f_n(x))|\ >\ 1\ \ \}.$$

Observe that $S(f_n)$ is a sub-complex of $Y$ and $S(f_{n+1})\subset S(f_n).$  We
claim if $S(f_n)\ne\phi$ then the cover $p_{n}$ may be chosen such that
$S(f_{n+1})$ is a proper subcomplex of $S(f_n).$ Since $Y$ is a finite complex it
then follows that for some
$n\ge0$ that
$S(f_n)$ is empty. Then $f_n$ is a $\pi_1$-injective embedding of $Y$ into the
finite cover $\tilde{M}_{n}$ which proves the theorem.

To prove the claim, suppose that $a,b\in Y$ are distinct points and $f_n(a)=f_n(b).$ Let $\gamma$ be
a path in $Y$ from $a$ to $b.$ Then $\alpha = [f_n\circ\gamma]\in \pi_1(\tilde{M}_n)\le\pi_1(M),$
and since $Y$ is a subspace of a covering of $M$ it is clear that $\alpha$ is
non-trivial. Since $G$ is separable there is a finite index subgroup $H <
\pi_1(M)$ which contains $G$ but does not contain $\alpha.$ The subgroup $H\cap
\pi_1\tilde{M}_{n}$ has finite index in $\pi_1(\tilde{M}_n).$ Let
$p_{n+1}:\tilde{M}_{n+1}\rightarrow \tilde{M}_n$ be the cover corresponding to this subgroup.
  It is
clear that $f_n:Y\rightarrow\tilde{M}_n$ lifts to $f_{n+1}:Y\rightarrow
\tilde{M}_{n+1}$ and, since $\alpha$ does not lift, that $f_{n+1}(a)\ne f_{n+1}(b).$
Thus
$S(f_{n+1})$ is a proper subsest of $S(f_n).$ This proves the claim.\qed

\medskip
From this one recovers the following well-known result:
\begin{corollary}[separable surface subgroups] Suppose that $f:S\rightarrow M$ is a continuous map of a closed surface with $\chi(S)\le0$ into an irreducible 3-manifold $M$ and suppose that $f_*:\pi_1S\rightarrow\pi_1M$ is injective. Suppose that $f_*(\pi_1S)$ is a separable subgroup of $\pi_1M.$ Then there is a finite cover $p:\tilde{M}\rightarrow M$ and an embedding $g:S\rightarrow \tilde{M}$ such that $p\circ g$ is homotopic to $f.$\end{corollary}
\demo By the theorem $f_*(\pi_1S)$ is virtually embedded so there is a finite cover $p:\tilde{M}\rightarrow M$ and an incompressible submanifold $Y\subset \tilde{M}$ such that $p_*(\pi_1Y)=f_*(\pi_1S).$ It is easy to see that we may choose $Y$ to be an embedding of $S$ into $\tilde{M}.$ Since $M$ is irreducible it is a $K(\pi_1M,1),$ hence the maps  $p|Y\equiv S\rightarrow M$ and $f:S\rightarrow M$ are homotopic.\qed

\begin{corollary} Suppose that  $M$ is a $3$-manifold which is homotopy equivalent to a tubed surface. Then every finitely generated subgroup of $\pi_1M$ is virtually embedded.
\end{corollary}

\begin{theorem}[totally geodesic is separable]\label{totgeodesic} Suppose that $M$ is a totally geodesic hyperbolic $k$-manifold immersed in a convex hyperbolic $n$-manifold $N$ with $k<n.$ Also suppose that $\pi_1N$ is finitely generated. Then $\pi_1M$ is a separable subgroup of $\pi_1N.$\end{theorem}
\demo This can be proved by an extension of  the method Long used in the case $k=2,n=3$ see \cite{Lo}.\qed

\section{Surfaces in 3-Manifolds.}

\medskip
We start with two results that will be needed for some of the $3$-manifold applications.
The following is well known, see for example theorem (3.15) of \cite{MT}. In dimensions $4$ and higher there are examples of two geometrically finite groups whose intersection is not finitely generated, see \cite{suss}.

\begin{theorem}[intersections of GF are GF]\label{intersectgeomfinite}  If $\Gamma_1$ and $\Gamma_2$ are geometrically finite subgroups of a discrete subgroup of $Isom({\mathbb H}^3)$ then $\Gamma_1\cap\Gamma_2$ is geometrically finite.\end{theorem}

A group $G$ has the {\it finitely generated intersection
property} or {\it FGIP} if the intersection of two finitely generated subgroups
is always finitely generated. The following result in this form is due to Susskind \cite{su};
see also Hempel  \cite{HE} and \cite{MT} corollary (3.16).

\begin{theorem}[GF infinite covolume implies FGIP]\label{FGIP} If $G$ is a geometrically finite Kleinian group of infinite co-volume then $G$ has the FGIP.\end{theorem}

The fundamental group of a hyperbolic surface bundle does not have FGIP, so the hypothesis of infinite co-volume is necessary.

 \medskip
The hypothesis of the convex combination theorem that one can thicken without bumping means that there are restrictions on the cusps of the manifolds to be glued. In particular two rank-1 cusps in the same rank-2 cusp of a $3$-manifold cannot be thickened without bumping unless they are parallel.

Here is an algebraic viewpoint. Rank-1 cusps give ${\mathbb Z}$ subgroups, and two non-parallel cusps in the same rank-2 cusp will generate a ${\mathbb Z}\oplus{\mathbb Z}$ group in the group, $G,$  generated by the two subgroups. Thus $G$ is not an amalgamated free product of the subgroups.

 Suppose $M_1$ and $M_2$ are $3$-manifolds with rank-1 cusps which intersect in $M.$ We may first glue rank-2 cusps onto each rank-1 cusp of  $M_1$ and $M_2$.
  
\centerline{\epsfysize=40mm
\epsfbox{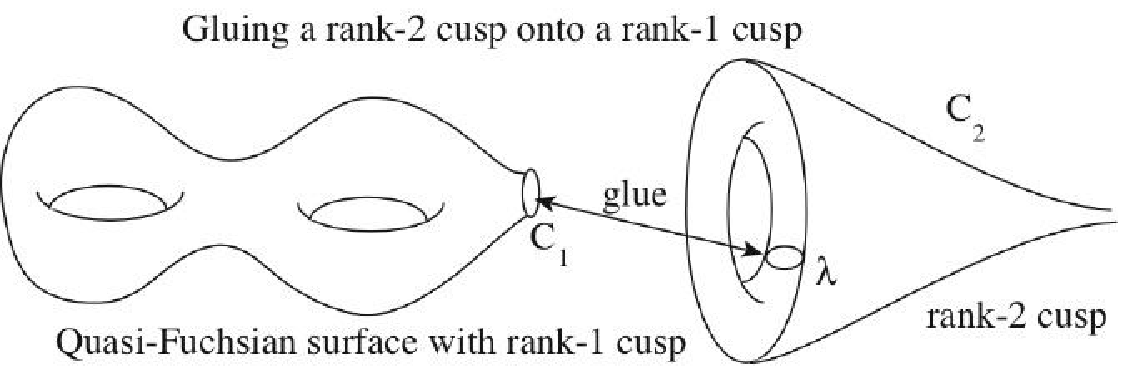}}   

  This corresponds to an amalgamated free product of a ${\mathbb Z}\oplus{\mathbb Z}$ with $\pi_1M_i.$ This produces two new geometrically finite manifolds $M_1^+,M_2^+$ which may now be glued.  The process of gluing a rank-2 cusp onto a rank-1 cusp can be done with  the convex combination theorem:

 \medskip
 \begin{theorem}[adding a rank-2 cusp]\label{addrank2cusp}
Suppose that $M$ is a convex hyperbolic $3$-manifold and that $f:N\rightarrow M$ is a locally-isometric immersion of a geometrically finite hyperbolic $3$-manifold $N.$ Suppose that $C_2$ is a rank-2 cusp of $M$ and   $C_1$ is a component of $f^{-1}(C_2)$ which is a rank-1 cusp of $N.$  Then there is a finite cover $\tilde{C}_2$ of $C_2$ and a geometrically finite $3$-manifold $N^+ = N \cup \tilde{C}_2$ with $N\cap\tilde{C}_2=C_1,$ where $C_1\subset N$ is identified with a subset of $\tilde{C}_2$ using a lift of $f.$   
\end{theorem}
\demo This follows from (\ref{increaserank}).\qed
 
 \medskip{\bf Remark.} It is easy to extend the above result to the setting where one has finitely many convex manifolds and finitely many rank-2 cusps, and it is required to glue some of the rank-1 cusps to cyclic covers of the rank-2 cusps in such a way that if more than one rank-1 cusp is glued onto the same rank-2 cusp then the rank-1 cusps are glued along parallel curves in the boundary of the rank-2 cusps.

\begin{corollary} [GF  tubed surfaces]\label{tubedsurface} Suppose $M$ is a compact  $3$-manifold whose interior admits a complete hyperbolic metric.  Suppose that $S$ is a compact, connected surface with $\chi(S)<0$ and $f:S\rightarrow M$ is $\pi_1$-injective such that each component of $\partial S$ is mapped into some torus boundary component of $M.$ Suppose that $f_*(\pi_1S)$ is a geometrically-finite subgroup of $\pi_1M.$  Let
$S^+$ denote the tubed surface obtained by gluing one torus onto each boundary
component of $S.$ Then $f$ extends to a $\pi_1$-injective map $f:S^+\rightarrow M$ such that
such that $f_*(\pi_1S^+)$ is a geometrically-finite subgroup of $\pi_1M.$
\end{corollary}
\demo We apply the previous theorem to add rank-2 cusps to the boundary components of $S$ one at a time.\qed

\medskip
{\bf Definition.} Suppose that $N$ is a convex hyperbolic 3-manifold of finite volume which equals its convex core and that $N$ has spine a tubed surface $S^+.$  In what follows we will the use of the term {\em tubed surface} to refer to either the 2-complex $S^+,$ or  to the geometrically finite hyperbolic 3-manifold $N$  and write this as $Core(S^+).$ 

\medskip
Suppose that  $\tilde{N}\rightarrow N$ is a finite cyclic cover to which $S$ lifts, then $\tilde{N}$ contains a closed surface, $2S,$ homeomorphic to the double of $S$ along its boundary. If the cover has degree bigger than $1$ then this surface may be chosen such that it is not homotopic into the boundary of a compact core of $\tilde{N}.$   Clearly $2S$ is $\pi_1$-injective in $\tilde{N}.$  In particular (\ref{tubedsurface}) implies that $\pi_1M$ contains a non-perpiheral surface group.

\medskip
Let $N=Core(S^+).$ In \cite{CL1} it  was observed that given a rank-2 cusp, $C,$ of $N$ every sufficiently large Dehn-filling of $C$ can be given a Riemannian metric of negative sectional curvature which is hyperbolic outside of $C.$ Suppose $f:N\rightarrow M$ is a local isometry and $N^+$ is a Dehn-filling of $N$ along $C.$ The same filling done on the corresponding cusp of $M$ then gives a local isometry $N^+\rightarrow M^+$ of Dehn-filled manifolds. If this is done to all the cusps of $N$ then, since $N$ is convex, this map is $\pi_1$-injective. This was the main step in the proof of the main theorem of \cite{CL1}.  It also leads to a quick proof of a virtual-Haken Dehn-filling result of the type in \cite{CL2}.

\medskip
The following is the main tool we need for studying immersed boundary slopes.

\begin{theorem}[gluing tubed surfaces]\label{tubeglue}  Suppose that $W$ is a convex hyperbolic 3-manifold and for $i\in\{1,2\}$ that $M_i$ is a tubed surface and $f_i:(M_i,m_i)\rightarrow (W,w_0)$ is a local isometry. Set $\Gamma_i=\pi_1(M_i,m_i).$ Then there are finite covers $p_i:\tilde{M}_i\rightarrow M_i$ and a simple gluing $\tilde{M}=S(f_1\circ p_1,f_2\circ p_2)$ of $\tilde{M}_1,\tilde{M}_2$ such that $\tilde{M}$ has a convex thickening.   Furthermore $\Gamma_i^-=p_{i*}(\pi_1(\tilde{M}_i,\tilde{m}_i))$ contains $\Gamma=f_{1*}(\Gamma_1)\cap f_{2*}(\Gamma_2).$ Also $\pi_1\tilde{M}$ is a free product of $\Gamma_1^-$ and $\Gamma_2^-$ amalgamated along $\Gamma.$
\end{theorem}
\demo  Since tubed surfaces are geometrically finite then, by theorem (\ref{intersectgeomfinite}), the group $\Gamma$ is geometrically finite  and thus finitely generated.  By  lemma (\ref{tubedisLERF}) $\Gamma_1$ and $\Gamma_2$ are LERF.  Thus $\Gamma$ is separable in both $\Gamma_1$ and $\Gamma_2.$ Define $N_i=T_{\kappa}(M_i)$ then $f_i$ extends to a local isometry $f_i:N_i\rightarrow W.$  Applying the virtual simple gluing theorem  to $N_1,N_2$ mapped into $W,$ it follows that there are finite covers $p_i:Y_i\rightarrow N_i$  which have a simple gluing $Y.$ Let $p_i|:\tilde{M}_i\rightarrow M_i$ be the restriction of the covering $p_i.$
We can now apply the convex combination theorem to $\tilde{M}_1\subset Y_1$ and $\tilde{M}_2\subset Y_2$ and deduce that $\tilde{M}$ has a convex thickening. 
\qed

\medskip
{\bf Definition.} A compact,  connected, orientable $3$-manifold, $M,$ is a {\it book of $I$-bundles} if it contains a submanifold $N$ such that each component of $N$ is a solid torus and the closure, $W,$ of $M\setminus N$ is an $I$-bundle over a compact (not necessarily connected) surface, $F,$ which contains no component $F_0$ with $\chi(F_0)>0.$ We also assume that $W\cap N$  is the sub $I$-bundle over $\partial F$ and is $\pi_1$-injective in $N$ (and therefore in $M$). It easily follows that $M$ is irreducible and boundary-irreducible and that every sub-bundle of $W$ is $\pi_1$-injective in $M.$  A simple-gluing of two surfaces along an incompressible subsurface gives a two-complex which is the spine of a book of $I$-bundles.

In the case there are no parabolics the next result is similar to a special case of corollary 5 of Gitik's paper \cite{G3}, and also (with a little work) to the combination theorem of Bestvina-Feighn \cite{BF}.  

\begin{corollary}[virtual surface gluing]\label{GFsurfaceglue} Let $S_1$ and $S_2$ be geometrically finite surfaces in a convex hyperbolic 3-manifold $M.$ Suppose that for each rank-2 cusp $C$ of $M,$  every rank-1 cusp of $S_1$ in $C$ is parallel to every rank-1 cusp of $S_2$ in $C.$ Then for all compatible choices of basepoints, there are finite covers of the convex hulls of $S_1$ and $S_2$ that have a simple gluing, $N,$ which is a book of $I$-bundles and is a geometrically finite subgroup. 
\end{corollary}
\demo Create tubed surfaces $Core(S_1^+),Core(S_2^+),$ from  $S_1$ and $S_2.$   Glue finite covers of these using (\ref{tubeglue}) to obtain a geometrically finite 3-manifold $P.$ Using the hypothesis on parallel cusps one can throw away the rank-2 cusps in $P.$ This last step can be accomplished by taking a covering of $P$ corresponding to a subgroup that only contain  parabolics parallel to the cusps of $S_1$ and $S_2.$\qed

\begin{corollary}\label{corol1} Suppose that $M$ is a closed hyperbolic 3-manifold and $S$ is a
closed, connected surface with $\chi(S)<0.$ Suppose that
$f:S\rightarrow M$ is $\pi_1$-injective and not homotopic to an embedding. Also
suppose that $f_*(\pi_1S)$ is a maximal surface subgroup of $\pi_1M.$ Then there is a book of $I$-bundles, $N,$ with $\beta_2(N)\ge2$ which isometrically immerses into $M.$ Thus for all $n>0$ there is a subgroup $G$ of finite index in $\pi_1N$
which is a geometrically finite subgroup of $ \pi_1M$ and
such that
$\beta_2(G) > n$ and $G$ is freely indecomposable.
\end{corollary}
\demo The surface $S$ is either a virtual fiber or quasi-Fuchsian. We first consider the case
 that $S$ is quasi-Fuchsian.  Since $S$ is not homotopic to an embedding there are two conjugates, $A\ne B$ of $f_*(\pi_1S)$ such that $C=A\cap B\ne 1.$ Given a Kleinian group $H$ we denote its limit set by $\Lambda(H).$ Since $A$ is quasi-Fuchsian $\Lambda(A)\cong S^1.$ The subgroup of $\pi_1M$ which stabilizes $\Lambda(A)$ is a surface group. By maximallity this group is $A.$ Since $A\ne B$ it follows that $\Lambda(A)\ne\Lambda(B).$ By \cite{sw}
 $\Lambda(A\cap B)= \Lambda(A)\cap\Lambda(B).$ This is therefore a proper non-empty subset of $\Lambda(A)$ and therefore not homeomorphic to a circle. Hence $C = A\cap B$ is not the fundamental group of a closed surface. Thus $C$ has infinite index in $A$ and in $B.$
 
  We now apply theorem (\ref{GFsurfaceglue}) to $A$ and $B.$ Thus there are finite covers $S_1,S_2$ of $S$ corresponding to the subgroups $A',B'$ and $\pi_1$-injective proper subsurfaces $T_i\subset S_i$ corresponding to $C.$ Let $X$ be the 2-complex obtained by gluing the two closed surfaces $S_1$ to $S_2$ by identifying $T_1$ with $T_2$ then $X$ is a topological realization of $A'*_CB'.$ Furthermore $X$ is the spine of a book of $I$-bundles.
 
There is a $n$-fold cyclic cover $\tilde{X}\rightarrow X$ for which $S_1$ lifts.  The $n$ distinct lifts of $S_1$ and the pre-image of $S_2$ show $\beta_2(\tilde{X})\ge n+1.$ Since $\tilde{X}$ is also the spine of a book of $I$-bundles it is clear that $\pi_1\tilde{X}$ is freely indecomposable.

The remaining case is that $f:S\rightarrow M$ is a virtual fiber. Thus there is a finite cover
$\tilde{M}\rightarrow M$ and a  lift
$\tilde{f}:S\rightarrow
\tilde{M}$ which is homotopic to an embedding. Furthermore the image is a fiber
of a fibration of $\tilde{M}$ over the circle. After taking a further finite cover  we may perform the cut and cross-join construction along a non-separating curve in $S$ (see \cite{CLR2})  to produce a $\pi_1$-injective  surface $S'$ immersed in $\tilde{M}$ with $\chi(S') = \chi(S)$ and $[S]=[S']\in H_2(\tilde{M}).$
 
Since $S'$ is transverse to the suspension flow of the fibering it follows that $S'$ is not homotopic to an embedding in $\tilde{M'}.$ If $S'$ is quasi-Fuchsian then it is a maximal surface subgroup of $\pi_1\tilde{M}$ because $[S']=[S]\in H_2(\tilde{M}).$ But the fiber $S$ has the
smallest norm in its homology class. If the surface subgroup $S'$ were not maximal there would be a surface $S''$ which is immersed into $\tilde{M}$ and then $[S]=k\cdot [S'']$ for some $k>1.$ Since $S$ is connected this is a contradiction. Hence $S'$ gives a maximal surface group in $\pi_1\tilde{M}.$ The result now follows from the first case. 

 There remains the possibility that $S'$ is also a virtual fiber. In this case there is a finite cover of $\tilde{M}$ in which $S'$ lifts to an embedding and this cover increases the rank of the first homology of $M.$ We can repeat this procedure either producing a maximal quasi-Fuchsian surface group in some finite cover which is not homotopic to an embedding (and obtaining the desired result) or else increasing the rank of the first homology of $M$ by an arbitrary amount.
 
 Suppose that $S'$ is a compact, connected,  closed, oriented  surface embedded in a 3-manifold and of minimal Thurston norm in its homology class. Also suppose that the homology class of $S'$ is not in the cone on the interior of some top dimensional face of the unit ball of the Thurston-norm then $S'$ is not a fiber and is thus quasi-Fuchsian. It follows that we can find arbitrarily many such surfaces in distinct homology classes in some finite cover of $M.$ Projecting these into $M$ gives a collection of immersed surfaces, which are possibly no longer maximal surface groups in $M.$ Each such surface has finite index in some other maximal surface group. These subgroups are therefore quasi-Fuchsian.  The corresponding surfaces can not all be homotoped to be embedded and pairwise disjoint by the Haken finiteness theorem.  This gives either one  maximal quasi-Fuchian surface which cant be homotoped to be embedded or two quasi-Fuchsian surfaces $S_1,S_2$ which cannot be homotoped to be disjoint. The result follows as before.
 \qed

\begin{theorem}[gluing two rank-1 cusps: algebraic version]
Suppose that $G$ is a torsion-free Kleinian group. Suppose that $H$ is a geometrically finite subgroup of $G$ and that $P,P'$ are two maximal parabolic subgroups in $H$ each of which is infinite cyclic.  Suppose that $T$ is  a maximal rank-2 parabolic subgroup of $G$ which contains  $P$ and there is $\gamma\in G$ such that $P'=\gamma P\gamma^{-1}.$ Then there is $n>0$ with the following property. Suppose  that $\alpha\in T$ and the subgroup of $T$ generated by $P$ and $\alpha$ has finite index at least $n.$ Set $\beta=\gamma\alpha$ then the subgroup of $G$ generated by $H$ and $\beta$ is the HNN extension $H*_P$ amalgamated along $P$ and $P'$ given by $< H,\beta :\ \beta P\beta^{-1} = P'>.$  Futhermore this group is geometrically finite.\
\end{theorem}

\medskip
A geometric formulation of this result involves a generalization of the notion of {\em spinning}  an annulus boundary component of a surface  around a torus boundary component of a 3-manifold that contains the surface. This idea was introduced in \cite{FF} and used in \cite{CL1}, \cite{CL2}. We recall the construction. Suppose $S$ is a compact surface immersed in a 3-manifold $M$ and that $\partial S$ has two components, $\alpha$ and $\beta,$ which both lie on a torus $T\subset \partial M.$ We suppose the immersion maps $\alpha$ and $\beta$ to the same loop on $T.$ Attach an annulus $A$ to the boundary of $S$ and choose an immersion of the annulus into $M$ so that it wraps some number $n\ge0$ times around $T.$ We describe this process by saying the two boundary components $\alpha$ and $\beta$ have been glued using an annulus that spins $n$ times around the torus $T.$
 
 \centerline{\epsfysize=60mm
\epsfbox{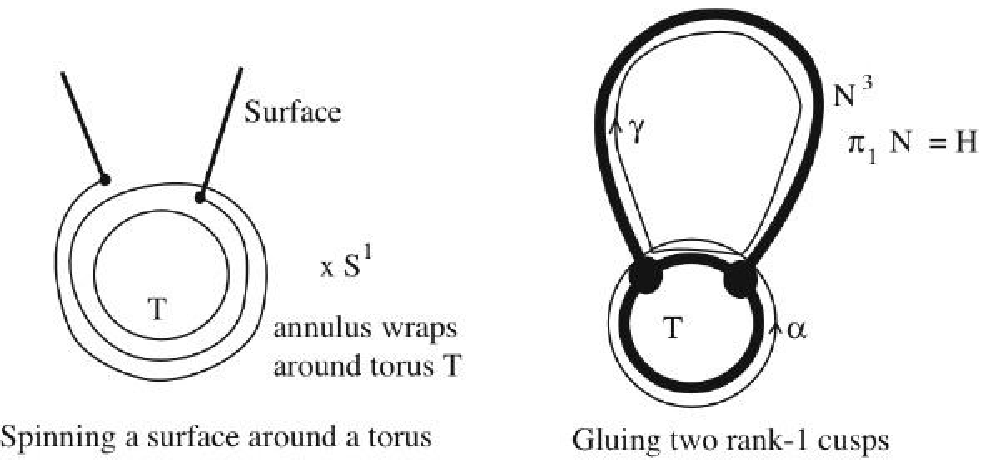}}      

This operation can be extended to the case where we are given an immersion $f:N\rightarrow M$ of a 3-manifold $N$  into $M$ and there are two annuli $A_0,A_1\subset \partial N$ which both map $m$ times around  the same annulus in some torus $T\subset \partial M.$ To be more precise, if $f:N\rightarrow M$ is the immersion we require there is a homeomorphism $h:A_1\rightarrow A_0$ such that $f|A_1=(f|A_0)\circ h.$ Then we glue $A\times[0,1]$ to $N$ to obtain a 3-manifold $N^+$ by identifying $A\times i$ with $A_i$ for $i=0,1$ and choose an extension of the immersion $f$ over $A\times[0,1]$ which spins $A$  round $T$ a total of $n$ times. The meaning of this last statement is the following. Identify $A_0$ with $S^1\times[0,s]$ and $T$ with $S^1\times S^1$ in such a way that $f|A_0$ is given by the map $f(e^{i\theta},x)=(e^{i m\theta},e^{ix}).$  Then define $f:A\times[0,1]\rightarrow S^1\times S^1$ by $f((e^{i\theta},x),y)=(e^{i\theta},e^{i(x+2n\pi y)}).$ We describe $N^+$ and the resulting immersion $f:N^+\rightarrow M$ as obtained from the immersion of $N$ into $M$ by {\em gluing two annuli in $\partial N$ together after spinning $n$ times around the torus $T.$} We can now state a geometric reformulation of the preceding result:

\begin{theorem}[gluing two rank-1 cusps: geometric version]\label{spin} Suppose that $M$ is a compact 3-manifold whose interior admits a complete hyperbolic metric. Suppose that $N$ is a compact 3-manifold corresponding to a geometrically finite subgroup of $\pi_1M.$  Let $f:N\rightarrow M$ be the immersion corresponding to this.  Suppose that $\partial N$ contains two annuli which map the same number of times around an annulus  in a torus $T\subset\partial M.$ Then for $n>0$ sufficiently large, the result of gluing the two annuli after spinning $n$ times around $T$ gives an isometric immersion of a geometrically finite $3$-manifold $N^+$ into $M.$\end{theorem}
\demo Glue the two rank-1 cusps of $N$ onto a large cyclic cover of a small rank-2 cusp in $M.$ Now take the cover which unwraps the rank-2 cusps to give the desired 3-manifold.\qed

\medskip
In particular if $N$ is a quasi-Fuchsian 3-manifold then $N^+$ is just a  thickened version (regular neighborhood) of the surface obtained by spinning two boundary components of the surface round $T.$ The manifold $N^+$ contains an accidental parabolic.

\section{Multiple Immersed Boundary Slopes.}
{\bf Definition.} An {\em immersed slope} $\alpha$ on a torus $T$ is an element of $(H_1(T;{\mathbb Z})-0)/(\pm1).$ Thus $\alpha$ is uniquely represented by the homotopy class of an
unoriented loop, also denoted $\alpha,$ in $T$ which is not contractible. If $\alpha$ is homotopic in $T$ to a simple closed curve  then $\alpha$ is called a {\it  slope.}

 Suppose $M$ is a compact 3-manifold and
$T$ is a torus boundary component of $M$. 
An {\it immersed boundary
slope}  is an immersed slope $\alpha$ on $T$ such that there is a compact
surface $S$ with non-empty boundary and a
$\pi_1$-injective map $f:S\rightarrow M$ such that:
\nb for every boundary component $\beta\subset\partial S$  we have $f_*([\beta]) 
= [\alpha].$ 
\nb $f$ is not homotopic rel $\partial S$ to a map into $\partial M.$ 

\medskip  If $S$ is embedded we say that
$\alpha$ is a {\em boundary slope.} If
$S$ does not lift to a fiber of a fibration of a finite cover $\tilde{M}$ then
$\alpha$ is a {\em strict boundary slope.} If $\alpha$ is a boundary slope which
is not a strict boundary slope and $S$ is embedded then either $S$ is  a fiber of a fibration of $M$
or else $S$ separates $M$ into two components each of which is a twisted
$I$-bundle over a surface. In this case we say that $S$ is a {\em semi-fiber.} If $\alpha$ is a
{\it  slope}  and for some $n>0$ the immersed slope $n\cdot\alpha$ is an immersed boundary slope we say that $\alpha$ is a {\it multiple immersed boundary slope} or MIBS.

\medskip Hatcher showed \cite{HA} that if $M$ is compact and has boundary a torus then there are only finitely many boundary slopes. 
In \cite{BC} it was shown that if $M$ is Seifert fibered then every immersed boundary slope is also a boundary slope and there are only two boundary slopes; the longitude (rationally null-homologous slope) and the slope of a fiber. The same result holds for MIBS. In this section we will show that if the interior of $M$ admits a complete hyperbolic metric then every slope is a MIBS. 

\medskip
Here is an outline of what follows. Suppose $S,S'$ are two quasi-Fuchsian surfaces each with two boundary components which correspond to different boundary slopes $\alpha,\beta.$ One would like to cut and cross join these surfaces along an arc connecting the two boundary components in each surface to produce a surface with boundary components $\alpha+\beta$ and $\alpha-\beta.$ Finally one would like to tube two copies of this surface together using the $\alpha-\beta$ slopes to obtain an immersed surface with two boundary components both with slope $\alpha+\beta.$ 

There are many problems with trying to do this directly, the most obvious one being that in general one can't expect to produce a $\pi_1$-injective surface this way. Instead we first add rank-2 cusps to the surfaces. Then finite covers of these tubed surfaces can be glued to give a geometrically finite manifold with a torus boundary component having the property that every slope on this torus is homologous to a cycle on the union of the other torus boundary components. This 2-chain is represented by an incompressible surface. Finally two copies of this surface are glued by spinning and gluing all the boundary components except those on the chosen torus. This gives the desired immersion.

\medskip
Culler and Shalen \cite{CS},\cite{CGLS} used character varieties to show:

\begin{theorem} Suppose that $M$ is a compact, orientable 3-manifold with boundary a torus $T$ and
that the interior of $M$ admits a complete hyperbolic structure of finite volume.  Then there are two strict boundary slopes on
$T.$\end{theorem}

Combining this with lemma (2.3)  of \cite{CL1} yields the following:
\begin{addendum}\label{2QFslopes}
There are two incompressible, $\partial$-incompressible, quasifuchsian surfaces
$S_1,S_2$ in $M$ with boundary slopes $\alpha_1\ne\alpha_2.$\end{addendum}

\medskip
The existence of the following cover is perhaps independently interesting.

\begin{theorem} \label{varyslope} Suppose that $M$ is a compact, orientable 3-manifold with boundary a torus $T$ and
that the interior of $M$ admits a complete hyperbolic structure of finite volume.  Suppose that
$\partial M = T$ is a torus. Then there is an infinite cover $p:\tilde{M}\rightarrow M$ such that $\pi_1\tilde{M}$ is finitely generated and there are distinct tori $T_1,\cdots,T_n\subset \partial\tilde{M}$ with the property that $$0=incl_*:H_1(T_1)\rightarrow H_1(\tilde{M},\cup_{i=2}^n T_i).$$
\end{theorem}

\demo By addendum (\ref{2QFslopes}) there are two quasi-Fuchsian surfaces $S_1,S_2$ in $M$ with distinct boundary slopes $\alpha_1,\alpha_2.$ By corollary (\ref{tubedsurface}) we may construct two geometrically finite manifolds $M_1,M_2$ with spines that are the tubed surfaces obtained from them.  Choose basepoints $m_i\in M_i$ which map to the same point in the cusp of $M.$  Then apply theorem (\ref{tubeglue}) to $\pi_1(M_1,m_1)$ and $\pi_1(M_2,m_2).$ This gives a geometrically finite manifold $\tilde{M}$ obtained by gluing finite covers of $M_1,M_2$ along a submanifold. 

There is a torus $T_1\subset\tilde{M}$ which corresponds to the rank-2 cusp  obtained by gluing the rank-2 cusps of $M_1$ and $M_2$ that contain the basepoints.
  Hence for $i\in\{1,2\}$ there is a component $\tilde{S}_i\subset \pi^{-1}S_i$  such that $T_1\cap\tilde{S}_i\ne \phi.$ Each component of  $T_1\cap\tilde{S}_i$ is a  loop $\tilde{\alpha}_i$ that covers $\alpha_i.$  By passing to a finite cover of $\tilde{M}$ we may assume $\tilde{\alpha}_1,\tilde{\alpha}_2$ generate $H_1(T_1).$ By passing to a further finite cover we may assume that $T_1\cap\tilde{S}_i$ is connected.  The result follows from consideration of the algebraic sum of $m$ copies of $\tilde{S}_1$ and $n$ copies of $\tilde{S}_2.$
 \qed
 
\begin{theorem}[All slopes are MIBS] Suppose that $M$ is a compact, orientable 3-manifold with boundary a torus $T$ and
that the interior of $M$ admits a complete hyperbolic structure of finite volume. Then there is a subgroup of finite index in $H_1(T;{\mathbb Z})$ such that every non-trivial element in this subgroup is an immersed boundary slope for a geometrically finite surface with exactly two boundary components. Thus every slope on $T$ is a MIBS. \end{theorem} 
\demo We apply (\ref{varyslope}) to obtain a cover $p:\tilde{M}\rightarrow M$ and a torus $T_1\subset\partial \tilde{M}$ with the property stated.  Let $n$ be the index of  $p_*H_1(T_1)$  in $H_1(T).$ Given an essential loop $\alpha$ on $T$ representing some slope, then $n\cdot\alpha$  lifts to a loop $\beta$ on $T_1.$ Thus there is a compact, connected, 2-sided, incompressible  surface $S$ properly embedded in $\tilde{M}$ such that  $S\cap T_1=\beta.$ For $i=0,1$ let $S_i$ be a copy of $S$ and $\beta_i\subset \partial S_i$ the boundary component corresponding to $\beta.$ For each boundary component $\gamma_0\subset \partial S_0$ with $\gamma_0\ne\beta_0$ attach the boundary components of an annulus to $\gamma_0$ and $\gamma_1$ to obtain a surface $R$ with two boundary components $\beta_0,\beta_1.$ Immerse $R$ into $\tilde{M}$ identifying the  two copies $S_i$ with $S$ and by spinning each annulus around the appropriate torus in $\partial\tilde{M}$ enough times to ensure the resulting immersed surface is $\pi_1$-injective and geometrically finite. That this can be done follows from (\ref{spin}). The composition $R\rightarrow \tilde{M}\rightarrow M$ is the desired surface.\qed

\gap

\noindent Mark Baker:\\Math Dept.\\Universit$\acute{e}$ de Rennes 1\\France\\email:
mark.baker@univ-rennes1.fr

\medskip
\noindent Daryl Cooper:\\Math Dept.\\UCSB\\Santa Barbara\\CA 93106\\USA\\email:
cooper@math.ucsb.edu


\begin{thebibliography}{99}

\bibitem{ALR} I. Agol, D.D. Long, A.W. Reid,
{\em The Bianchi groups are separable on geometrically finite subgroups,}
 Ann. of Math. (2) 153 (2001), no. 3, 599--621.
\bibitem{B} M.D. Baker,  {\em On boundary slopes of immersed incompressible surfaces,}
Ann. Inst. Fourier (Grenoble) 46 (1996), no. 5, 1443--1449.
\bibitem{BC} M.D. Baker, D. Cooper, {\em Immersed, virtually-embedded, boundary
slopes,} Topology Appl. 102 (2000), no. 3, 239--252.
\bibitem{BF} M. Bestvina, M. Feighn, {\em
A combination theorem for negatively curved groups,}
J. Differential Geom. 35 (1992), no. 1, 85--101.
\bibitem{Bow} B.H. Bowditch, {\em Geometrical finiteness for hyperbolic groups.} J. Funct. Anal. 113 (1993), no. 2, 245--317
{\em Addendum and correction to: "A combination theorem for negatively curved groups"} J. Differential Geom. 43 (1996), no. 4, 783--788.
\bibitem{CL1} D. Cooper, D.D. Long, {\em Some Surface Subgroups Survive
Surgery,}  Geom. Topol. 5 (2001), 347--367.
\bibitem{CL2} D. Cooper, D.D. Long, {\em
Virtually Haken Dehn-filling,}
J. Differential Geom. 52 (1999), no. 1, 173--187.
\bibitem{CLR2} D. Cooper, D.D. Long, A.W. Reid, {\em Bundles and finite
foliations,}  Invent. Math. 118 (1994), no. 2, 255--283.
\bibitem{CS} M. Culler, P.B. Shalen, {\em  
Boundary slopes of knots,}
Comment. Math. Helv. 74 (1999), no. 4, 530--547.
\bibitem{CGLS} M. Culler, C. McA. Gordon, J. Luecke, P.B.
Shalen, {\em  Dehn surgery on knots,} 
Ann. of Math. (2) 125 (1987), no. 2, 237--300.
\bibitem{FD} F. Dahmani, {\em Combination of convergence groups.} Geom. Topol. 7 (2003) 933-963
\bibitem{FF} B. Freedman, M.H. Freedman,  {\em Haken finiteness for
bounded $3$\--manifolds, locally free groups and cyclic covers,}  
Topology 37 (1998), no. 1, 133--147.
\bibitem{G1} R. Gitik, {\em Graphs and separability properties of groups,} J. Algebra 188 (1997), no. 1, 125--143.
\bibitem{G2} R. Gitik, {\em
On the combination theorem for negatively curved groups,}
Internat. J. Algebra Comput. 6 (1996), no. 6, 751--760.\\
{\em Corrected reprint of  ``On the combination theorem for negatively curved groups.''} Internat. J. Algebra Comput. 7 (1997), no. 2, 267--276.\\
{\em Errata: ``On the combination theorem for negatively curved groups"} Internat. J. Algebra Comput. 7 (1997), no. 2, 265.
\bibitem{G3} R. Gitik, {\em Ping-Pong on Negatively Curved Groups,} J. Algebra 217 (1999), 65--72
\bibitem{HRW} J. Hass, J.H. Rubinstein, S. Wang, {\em Boundary
slopes of immersed surfaces in 3-manifolds,} J. Differential
    Geom. 52 (1999), no. 2, 303--325. 
\bibitem{HWZ} J. Hass, S. Wang, Q. Zhou, {\em On finiteness of the
number of boundary slopes of immersed surfaces in
    3-manifolds,} Proc. Amer. Math. Soc. 130 (2002), no. 6, 1851--1857
\bibitem{HA} A. Hatcher, {\em On the boundary curves of incompressible surfaces,} Pacific J.
Math. 99 (1982) 373-377.
\bibitem{HE} J. Hempel, {\em The finitely generated intersection property
for Kleinian groups,}  Knot theory and manifolds (Vancouver, B.C., 1983),
18--24,  Lecture Notes in Math., 1144,  Springer, Berlin, 1985.
\bibitem{HE2} J. Hempel,  {\em 3-Manifolds,}
Ann. of Math. Studies, No. 86. 
Princeton University Press, Princeton, N. J.; University of Tokyo Press, Tokyo,
(1976).
\bibitem{KR} E. Kang, J.H. Rubinstein, {\em Ideal triangulations of 3-manifolds I: spun normal surface theory,} G \& T monographs, (2004), vol 7, 235--265
\bibitem{Lo} D.D. Long, {\em Immersions and embeddings of totally geodesic surfaces.} 
Bull. London Math. Soc. 19 (1987), no. 5, 481--484.
\bibitem{MA} J. Maher, {\em Virtually embedded boundary slopes,} Topology Appl. 95
(1999), no. 1, 63--74. 
\bibitem{Mal} A. I. Malcev, {\em On isomorphic matrix representations of infinite groups,} 
Mat. Sb. 8, 405-422 (1940).
\bibitem{MAS} B. Maskit, {\em  On Klein's combination theorem. IV, }
Trans. Amer. Math. Soc. 336 (1993), no. 1, 265--294.
\bibitem{MT} K. Matsuzaki, M. Taniguchi, {\em Hyperbolic manifolds and Kleinian Groups,}
Oxford Science pub. (1998)
\bibitem{OE} U. Oertel, {\em Boundaries of $\pi\sb 1$-injective surfaces,}
Topology Appl. 78 (1997), no. 3, 215--234. 
\bibitem{SC1} P. Scott, 
{\em Subgroups of surface groups are almost geometric,}
J. London Math. Soc. (2) 17 (1978), no. 3, 555--565.\\
{\em Correction to: ``Subgroups of surface groups are almost geometric''} 
J. London Math. Soc. (2) 32 (1985), no. 2, 217--220.
\bibitem{Serre} J.P. Serre, {\em Trees.} Springer-Verlag, Berlin, 2003
\bibitem{suss} P. Susskind, {\em An infinitely generated intersection of geometrically finite hyperbolic groups.} 
Proc. Amer. Math. Soc. 129 (2001), no. 9, 2643--2646
\bibitem{sw} P. Susskind, G.A. Swarup 
{\em Limit sets of geometrically finite hyperbolic groups.} 
Amer. J. Math. 114 (1992), no. 2, 233--250.
\bibitem{su} P. Susskind, {\em Kleinian groups with intersecting limit sets.}
J. Analyse Math. 52 (1989), 26--38.
\bibitem{Th1} W.P. Thurston, {\em Three-Dimensional Geometry and Topology,}
Princeton Univ. Press, (1997), Ed Silvio Levy
\bibitem{Th2} W.P. Thurston, {\em The Geometry and Topology of Three-Manifolds,}
Princeton notes, 1977. at http://www.msri.org/publications/books/gt3m/
\end{thebibliography}
\end{document}